\documentclass[final]{siamltex}

\pdfoutput=1

\usepackage{amsmath,amsfonts}
\usepackage{graphicx}
\usepackage{url}
\usepackage{subfig}

\newcommand{\mA}{\mathbf{A}}

\newcommand{\mJ}{\mathbf{J}}
\newcommand{\mK}{\mathbf{K}}
\newcommand{\mC}{\mathbf{C}}
\newcommand{\mX}{\mathbf{X}}
\newcommand{\mV}{\mathbf{V}}

\newcommand{\va}{\mathbf{a}}

\newcommand{\vc}{\mathbf{c}}
\newcommand{\vf}{\mathbf{f}}
\newcommand{\vj}{\mathbf{j}}
\newcommand{\vs}{\mathbf{s}}
\newcommand{\vt}{\mathbf{t}}
\newcommand{\vx}{\mathbf{x}}
\newcommand{\vu}{\mathbf{u}}
\newcommand{\vv}{\mathbf{v}}
\newcommand{\vy}{\mathbf{y}}
\newcommand{\vz}{\mathbf{z}}

\newcommand{\bmat}[1]{\begin{bmatrix}#1\end{bmatrix}}

\newcommand{\mini}[1]{\underset{#1}{\operatorname{minimize}\;}}

\newcommand{\st}{\text{subject to }}

\title{Input Subspace Detection for Dimension Reduction in High Dimensional Approximation}

\author{
Paul G.~Constantine\thanks{Stanford University, 
Stanford, California 94305 ({\tt paul.constantine@stanford.edu}).}
\and 
Qiqi Wang\thanks{Department of Aeronautics and Astronautics,
Massachusetts Institute of Technology, Cambridge, Massachusetts 02139 ({\tt qiqi@mit.edu}).}
}

\begin{document}

\maketitle

\begin{abstract}
\emph{This manuscript is superseded by Constantine, Dow, and Wang's ``Active Subspaces in Theory and Practice: Applications to Kriging Surfaces'' [SIAM J.~of Sci.~Comput., 36 (2014), pp.~A1500--A1524].}

Many multivariate functions encountered in practice vary primarily along a few directions in the space of input
parameters. When these directions correspond with coordinate directions, one may apply global sensitivity measures to
determine the parameters with the greatest contribution to the function's variability. However, these methods perform
poorly when the directions of variability are not aligned with the natural coordinates of the input space. We present a
method for detecting the directions of variability of a function using evaluations of its derivative with respect to the
input parameters. We demonstrate how to exploit these directions to construct a surrogate function that depends on
fewer variables than the original function, thus reducing the dimension of the original problem. We apply this procedure
to an exercise in uncertainty quantification using an elliptic PDE with a model for the coefficients that depends on 250
independent parameters. The dimension reduction procedure identifies a 5-dimensional subspace suitable for constructing
surrogates. 
\end{abstract}


\pagestyle{myheadings}
\thispagestyle{plain}
\markboth{P.~G. CONSTANTINE AND Q. WANG}{INPUT SUBSPACE DETECTION}

\section{Introduction \& Motivation}
In modern science and engineering practice, computational simulation is routinely employed to help test hypotheses
and explore new designs. As the speed and capability of computers increase, so does the complexity of simulations
through greater resolution and higher fidelity physical models. Expensive simulations requiring extensive time on
massive supercomputers are now commonplace. Due to the cost of these high-fidelity simulations, one often wishes to
approximate the output at many points in the space of inputs using a surrogate function or a meta-model. The parameters
of the surrogate are tuned with a budget-constrained number of costly high-fidelity runs, and the tuned surrogates are
used to study sensitivities or uncertainties in the simulation output with respect to variation in the input parameters.

However, many surrogate models suffer from the so-called curse of dimensionality. Loosely speaking, the work
required to construct and evaluate an accurate surrogate increases exponentially as the dimension of the parameter
space increases. For example, this curse limits the applicability of polynomial-based surrogates to problems with a
handful of input parameters. Even methods whose application is independent of the dimension of the parameter space --
such as radial basis functions or Gaussian process models -- often perform poorly if the function is not sufficiently
smooth and the training data are too sparse.

Fortunately, in many problems of interest with high dimensional input spaces, the output often depends on only a few
important parameters. Specifically, the variability in the output can be attributed to a subset of the inputs.
Surrogates can be adjusted to take advantage of this anisotropic parameter dependence. A common approach --
known as global sensitivity analysis~\cite{Saltelli08} -- involves a strategy for ranking the input variables and
biasing the choice of design points to capture the function's behavior as the important parameters are varied. In some
cases, the ranking procedure can use a priori knowledge from the mathematical model.
In other cases, it requires exploration of the output through sampling. In either case, methods based on
variance-based decompositions~\cite{Liu06} or high-dimensional model representations~\cite{Li02} choose a few important
parameters from the full set of inputs. For problems encountered in practice, this often
results in a \emph{dimension reduction} of the input space; a surrogate can be constructed on a function of fewer
variables with significantly less work.

In this paper, we present a generalization of subset selection methods. Namely, we seek a low-dimensional linear
subspace of the input parameter space that captures the majority of the output's variability. The subspace induces a
reduced set of coordinates, and surrogate functions can be trained on the reduced coordinates to approximate the output
in the full space. 

More precisely, for a function of interest $f=f(\vs)$ with $\vs\in\mathbb{R}^d$, we seek a function
$g=g(\vs_a)$ with $\vs_a\in\mathbb{R}^a$ such that $f\approx g$ with $a<d$. The approximate function takes the form
$g(\vs_a)=f(\mA\vs_a)$, where $\mA$ is a $d\times a$ matrix representing a linear map from $\mathbb{R}^a$ to
$\mathbb{R}^d$. Note that the introduction of $g$ is primarily for notation; each evaluation of $g$ is ultimately an
evaluation of $f$ at specially chosen input values. However, the dependence of $g$ on fewer variables makes it more
amenable to surrogate approximation.

A similar idea is proposed in~\cite{Lieberman10} in the
context of model reduction for inverse problems, but the method for computing the basis vectors that define the
subspace employs the residual of a system of equations representing a physical model. Our method applies to more general
multivariate functions, and it is particularly efficient if one can easily compute gradients of outputs with respect to
inputs. We discuss strategies for the case when only function evaluations are available, including an intriguing idea of
using new matrix completion techniques on a partially sampled matrix of finite difference approximations of the
gradient. 

\section{Input subspace detection and dimension reduction}
We assume that a given multivariate function $f(\vs)$ with $\vs\in\mathbb{R}^d$ varies primarily along a few
directions in the input space. However, these directions may not be aligned
with the natural coordinate system. The goal in this section is to construct a function $g$ that approximates 
$f$ but takes only as many inputs as directions of variability.
Our strategy is to first determine the directions along which $f$ varies most prominently;
we rotate our coordinate system according to these directions. We then define $g$ to depend on the subset of these
rotated coordinates that contain the majority of the variability in $f$. 

\subsection{Directions of variability}
Let $\Omega$ be a hyperrectangle defined by the vectors $\vs_l$ and $\vs_u$,
\begin{equation}
\Omega = \{\vs \;:\; \vs\in\mathbb{R}^d,\,\vs_l \leq \vs \leq \vs_u \}
\end{equation}
We assume without loss of generality that
$\mathbf{0}\in\mathbb{R}^d$ is the center of mass of $\Omega$. Define a scalar function $f:\Omega\rightarrow\mathbb{R}$
that takes $d$ inputs. Denote an element of $\Omega$ by a $d$-vector $\vs=(s_1,\dots,s_d)^T\in\Omega$. For the
analysis, we assume that $f$ is analytic in a region containing $\Omega$. Denote the $d$-vector $\vj=\vj(\vs)$ as
\begin{equation}
\vj^T \;=\; \nabla f \;=\; \bmat{\frac{\partial f}{\partial s_1}& \cdots & \frac{\partial f}{\partial s_d}},
\end{equation}
which is the Jacobian of $f$. Define the $d\times d$ matrix $\mC$ 
\begin{equation}
\label{eq:C}
\mC = \int_\Omega \vj \vj^T \,d\vs, 
\end{equation}
where we employ a shorthand $d\vs$ to denote a measure on $\Omega$. Note that $\mC$ is symmetric and positive
semidefinite, which implies it has an eigenvalue decomposition
\begin{equation}
\label{eq:eigdecomp}
\mC = \mV\Lambda \mV^T, \qquad \Lambda=\mathrm{diag}\,(\lambda_1,\dots,\lambda_d), \qquad \lambda_1\geq \cdots
\geq \lambda_d \geq 0.
\end{equation}
If $\vv_i$ is the $i$th column of $\mV$, then
\begin{equation}
\label{eq:eig}
\lambda_i \;=\; \vv_i^T \mC \vv_i 
  \;=\; \vv_i^T \left(\int_\Omega \vj \vj^T \,d\vs\right) \vv_i 
  \;=\; \int_\Omega \left(\vv_i^T \vj \vj^T \vv_i\right) \,d\vs 
  \;=\; \int_\Omega \left(\vj^T \vv_i\right)^2 \,d\vs.
\end{equation}
We examine the Taylor expansion of $f$ at the point $\vs+h\vv_i\in\Omega$ around the point $\vs$,  
\begin{equation}
\label{eq:taylor}
f(\vs+h\vv_i) - f(\vs) \;=\; h\,\vj(\vs)^T\vv_i + \dots.
\end{equation} 
Taking the root-mean-squared of \eqref{eq:taylor} and applying \eqref{eq:eig}, we get
\begin{equation}
\label{eq:bound}
\|f(\vs+h\vv_i)-f(\vs)\|_{L_2} \;=\; \mathcal{O}(h\sqrt{\lambda_i}),
\end{equation}
where $\|\cdot\|_{L_2}$ is the standard $L_2$ norm for functions defined on $\Omega$.
Note that \eqref{eq:bound} implies the following: if $\lambda_i=0$, then the function $f$ is constant along the
direction of $\vv_i$. We can use this flatness to construct a sampling strategy to approximate $f$ on a low
dimensional manifold of $\Omega$. 

\begin{figure}
\begin{center}
\includegraphics[scale=0.4]{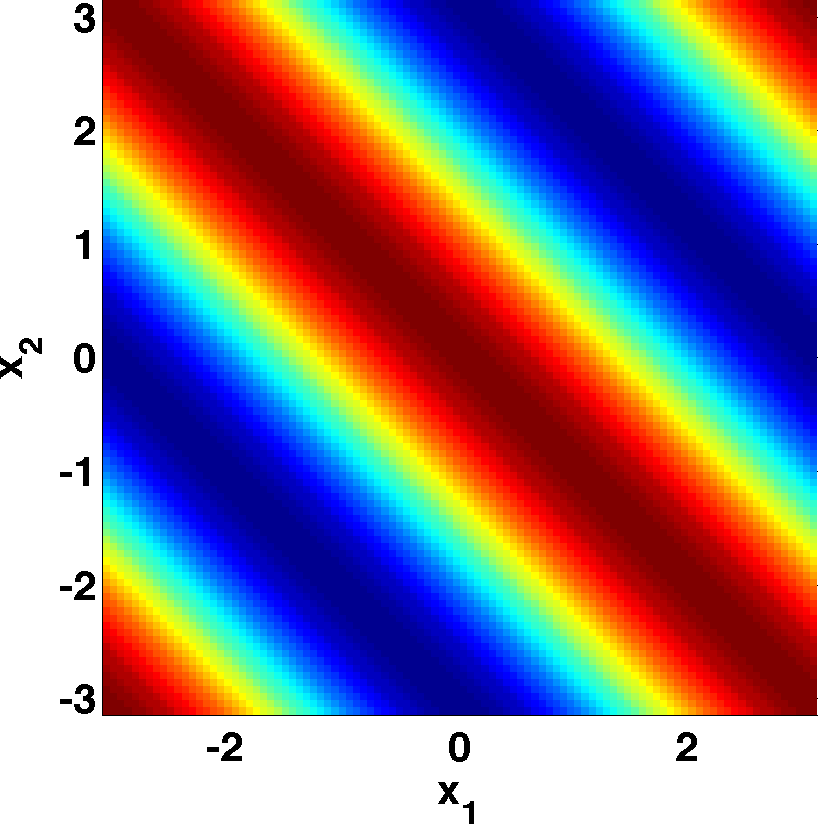}
\end{center}
\caption{The function $\cos(s_1+s_2)$ varies along the normalized vector $[\sqrt{2}/2, \sqrt{2}/2]^T$.}
\label{fig:cos}
\end{figure}

As an example, consider the function $f(\vs)=\cos(s_1+s_2)$ defined on $[-\pi,\pi]^2$; this function in plotted in
figure \ref{fig:cos}. The Jacobian of $f$ is
\begin{equation}
\nabla f = \bmat{-\sin(s_1+s_2) & -\sin(s_1+s_2)}.
\end{equation}
The matrix $\mC$ is then given by
\begin{equation}
\mC \;=\; \int_{-\pi}^\pi \int_{-\pi}^\pi 
\bmat{\sin^2(s_1+s_2) & \sin^2(s_1+s_2)\\ \sin^2(s_1+s_2) & \sin^2(s_1+s_2)}
\,ds_1 \,ds_2
\;=\;
2\pi^2\bmat{1 & 1\\ 1& 1}.
\end{equation}
The eigenvalue decomposition of $\mC$ from \eqref{eq:C} is
\begin{equation}
\mC \;=\;
\bmat{\sqrt{2}/2 & -\sqrt{2}/2 \\ \sqrt{2}/2 & \sqrt{2}/2}
\bmat{4\pi^2 & 0 \\ 0 & 0}
\bmat{\sqrt{2}/2 & \sqrt{2}/2 \\ -\sqrt{2}/2 & \sqrt{2}/2}. 
\end{equation}
Notice that the normalized vector $[\sqrt{2}/2, \sqrt{2}/2]^T$ -- the first eigenvector -- precisely identifies the
direction in the domain along which $f$ varies. Therefore, if we study $f$ along the line defined by $[\sqrt{2}/2,
\sqrt{2}/2]^T$, then we can understand the variation in $f$ over the whole domain through a projection.

\subsubsection{A note on ridge-type functions}
The previous function is an example of a \emph{ridge function}, which appear frequently in statistics~\cite{Cohen11}. A
ridge function takes the form
\begin{equation}
f \;=\; f(\va^T\vs) \;=\; f(t), \qquad \va\in\mathbb{R}^d, \quad \vs\in\Omega, \quad t=\va^T\vs.
\end{equation}
The Jacobian has a special form in this case:
\begin{equation}
\nabla f \;=\; \frac{df}{dt}\,\va^T.
\end{equation}
Then
\begin{equation}
\mC \;=\; \vv\lambda\vv^T, \qquad 
\lambda=\frac{\|\va\|^2}{a}\int \left(\frac{df}{dt}\right)^2\,dt,\qquad\|\vv\|=1,
\end{equation}
where the norm $\|\cdot\|$ is the
standard 2-norm on $\mathbb{R}^d$. The eigenvector $\vv$ is a normalized version of $\va$ that reveals the
direction of variability for $f$. In this case, $\mC$ has rank one when $df/dt\not=0$, and the vector $\vv$ can be
computed with one normalized point evaluation of $\vj(\vs)$.

\subsection{Dimension Reduction}

Assume that the eigenvalue decomposition \eqref{eq:eigdecomp} of $\mC$ can be partitioned as
\begin{equation}
\label{eq:eigpart}
\mV = \bmat{\mV_a & \mV_b},\qquad \Lambda = \bmat{\Lambda_a & \\ & 0},
\end{equation}
where $\mV_a$ has $a$ columns, $\mV_b$ has $b$ columns, and $a+b=d$. The columns of $\mV_a$ correspond to
directions along which $f$ varies, and the columns of $\mV_b$ correspond to directions along which $f$ is constant. We
can construct a rotated coordinate system since 
\begin{equation}
\label{eq:coord}
f(\vs) 
\;=\; f(\mV\mV^T\vs) 
\;=\; f(\mV_a\mV_a^T\vs + \mV_b\mV_b^T\vs) 
\;=\; f(\mV_a\vs_a + \mV_b\vs_b)
\;\equiv\; g(\vs_a,\vs_b). 
\end{equation}
By construction, the value of the function $g$ will not change as $\vs_b$ varies. One may be tempted to fix $\vs_b$
(say, set $\vs_b=0$) and treat $g$ as function of the $a$ variables $\vs_a$. However, there are two issues we must
address.

\subsubsection{Rotated coordinates}
First, what values can $\vs_a$ take? We can linearly transform the set of points $\vs\in\Omega$ to get a range for
$\vs_a$. In particular, we define the set 
\begin{equation}
\label{eq:reddom}
\Omega_a = \{\vs_a \;:\; \vs_a=\mV_a^T\vs,\; \vs\in\Omega\}. 
\end{equation}
Since $\Omega$ is convex, $\Omega_a$ is also convex, but this is about all we can say. The
coordinates $\vs_a$ cannot be varied independently within a set of independent intervals like a hyperrectangle, since
the transformed domain $\Omega_a$ will most likely not be a lower dimensional hypercube; imagine taking a photograph of
a rotated cube. For this reason, when we construct a surrogate on the reduced coordinates $\vs_a$, it must be flexible
enough to handle general convex domains in multiple dimensions; radial basis functions could be an appropriate choice.
We discuss sampling from the space $\Omega_a$ in section \ref{sec:linprog}.

\subsubsection{The domain of $f$}
We must ensure that all function evaluations of $f$ occur at points in the domain $\Omega$; each evaluation of
$g$ is ultimately an evaluation of $f$ at specially chosen input values. It is possible that the projection
$\mV_a\vs_a=\mV_a\mV_a^T\vs$ will not be in $\Omega$, and we do not want to assume anything about $f$ outside its
domain. Fortunately, we can take advantage of the flatness of $f$ to ensure that all evaluations occur within $\Omega$.
In short, for any point $\mV_a\vs_a$ that falls outside the domain of $f$, we can walk back along the directions in
which $f$ is constant until we reach a point in the domain. 

More precisely, if $\mV_a\vs_a\in\Omega$, then we evaluate $f$ at $\mV_a\vs_a$. If $\mV_a\vs_a\not\in\Omega$, then we
find $\vz\in\mathbb{R}^b$ such that $\mV_a\vs_a+\mV_b\vz\in\Omega$. Now we define $g$ as 
\begin{equation}
\label{eq:gfun}
g(\vs_a) \;=\; \left\{
\begin{array}{cl}
f(\mV_a\vs_a) & \mbox{ if $\mV_a\vs_a\in\Omega$, } \\
f(\mV_a\vs_a + \mV_b\vz) & \mbox{ if $\mV_a\vs_a\not\in\Omega$.} 
\end{array}
\right.
\end{equation}
Note that $\vz$ is often not uniquely determined, but we only need one for each deviant $\vs_a$. If $\Omega$ is a
hyperrectangle, then a $\vz$ can be found by solving a suitable linear program; see section \ref{sec:linprog}. We have
thus acheived our goal of constructing a function $g$ dependent on $a<d$ parameters that behaves like the $d$-variate
function $f$.

\begin{figure}
\begin{center}
\subfloat[$\cos(0.3s_1+0.7s_2)$]{
\includegraphics[scale=0.3]{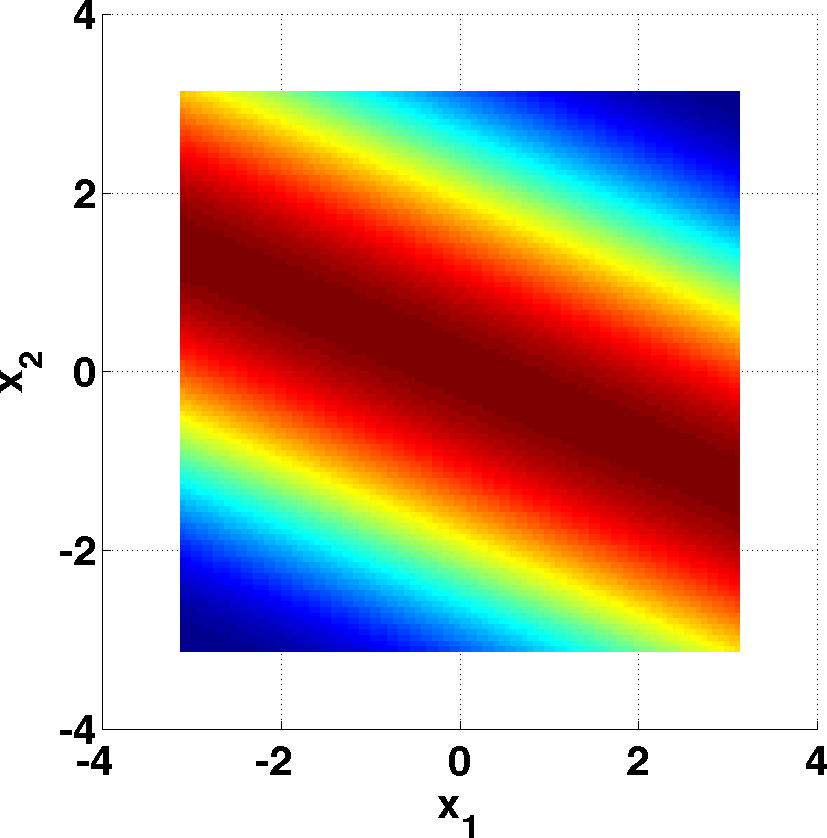}
}
\quad
\subfloat[Domains]{
\includegraphics[scale=0.3]{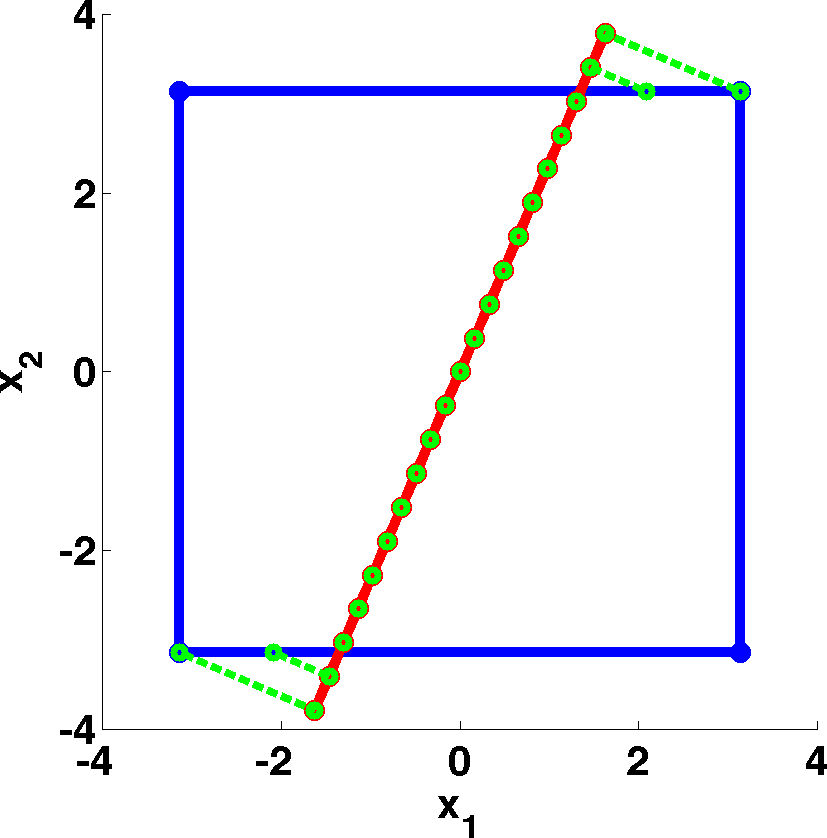}
}
\end{center}
\caption{The two-dimensional domain $[-\pi,\pi]^2$ of $\cos(0.3s_1+0.7s_2)$ shown in blue. The red shows the projection
of the domain onto the direction of variability. The green circles show the points where we evaluate $g$.}
\label{fig:domain}
\end{figure}

We demonstrate the rotation and reduction on a slight modification of the previous example. Let
$f(\vs)=\cos(0.3s_1+0.7s_2)$ be defined on $[-\pi,\pi]^2$ with gradient
\begin{equation}
\nabla f = \bmat{-0.3\sin(0.3s_1+0.7s_2) & -0.7\sin(0.3s_1+0.7s_2)}.
\end{equation}
Figure \ref{fig:domain} shows the domain $[-\pi,\pi]^2$ in blue. The projection of the domain onto the line
corresponding to the direction of variability of $f$ is shown in red. The red circles correspond to a possible sampling
of the reduced (one-dimensional) coordinates. Notice that some of the projected points fall outside $[-\pi,\pi]^2$ when
transformed back to the original two-dimensional space. The green circles show the points in $[-\pi,\pi]^2$ that are
substituted for the red points outside the domain when evaluating the function $g$. 

\section{Computational aspects}
We next consider four computational aspects for the dimension reduction procedure. We close this section with a
practical algorithm that summarizes the presentation.

\subsection{Low variability versus no variability}
In \eqref{eq:eigpart}, we assume that some of the eigenvalues are exactly zero. When this happens, the
estimate \eqref{eq:bound} tells us that $f$ is exactly constant along some directions. But what happens when the
eigenvalues are small but not zero? The estimate \eqref{eq:bound} addresses an averaged measure of variability
along a direction. Like any averaged measure, it does not preclude sharp, local variability. It is possible that one
could choose to ignore a direction because its associated eigenvalue is below a specified tolerance but subsequently
discover a sharp local feature in $f$ along this direction.

However, practical considerations like computational
budget often dominate the concerns when approximating functions in high dimensions. Any well-motivated strategy to
reduce cost is welcome. In this spirit, we treat the magnitudes of the eigenvalues as a ranking on the rotated
coordinates. If we desire an approximate bivariate function $g$ of $f$, then we choose the directions associated with
the two largest eigenvalues.

\subsection{Approximating the eigenvalues and eigenvectors}
For high dimensional functions found in practice, we expect that there will be a few dominant directions in the sense
described above. We usually are not able to compute the exact matrix $\mC$ from \eqref{eq:C}, but we can approximate
it. Assume for now that we can evaluate the exact Jacobian $\vj(\vs)^T$ given $\vs$. Then we can
approximate $\mC$ with a numerical quadrature rule. For simplicity, we use a Monte Carlo approximation.
For $i=1,\dots,k$, let $\vs_i$ be samples drawn from $\Omega$, and compute the $d\times k$ matrix
\begin{equation}
\label{eq:jacmat}
\mJ = \bmat{\vj(\vs_1) & \cdots & \vj(\vs_k)}.
\end{equation} 
Then 
\begin{equation}
\label{eq:ctilde}
\mC 
\;\approx\; 
\hat{\mC} 
\;=\; 
\frac{|\Omega|}{k} \sum_{i=1}^k \vj(\vs_i)\,\vj(\vs_i)^T 
\;=\; 
\frac{|\Omega|}{k}\,\mJ\mJ^T,
\end{equation}
where $|\Omega|$ is the volume of $\Omega$. The quality of the approximation can be controlled by the number of samples
$k$. For the Monte Carlo approximation, the variance of the approximation decreases like $k^{-1/2}$~\cite{Kalos08}.

Results from eigenvalue perturbation theory show that the error in the approximate eigenvalues is on
the order of the error in the matrix elements~\cite{GVL96}. More accurate numerical quadrature methods will result
in more accurate approximate eigenvalues, but many high order (e.g., interpolatory) multivariate quadrature rules
suffer from the same curse of dimensionality that we wish to avoid. For this reason, we rely on Monte Carlo methods. 

If $f$ is constant along some directions, then these directions will be in the null space of $\hat{\mC}$ when $k \geq
d$, i.e., when the number of Jacobian samples is greater than the number of parameters of $f$. The danger with the
approximation $\hat{\mC}$ is potentially overpredicting the dimension of the null space or, equivalently,
underpredicting the rank of $\mC$. In other words, the Jacobian evaluations at the design sites $\vs_i$ may indicate
that $f$ is flat along directions that it actually varies in $\Omega$. 

\subsection{Approximate Jacobians}
\label{sec:apjac}
Up to this point, we have assumed that the Jacobian $\vj(\vs)^T$ was available for computation. This is not true in many
cases, particularly if $f$ represents the output of a complex physical simulation. We therefore address the question of
approximating the Jacobian from point evaluations of $f$. 

If $f$ can be evaluated at will, then a finite difference approximation along the original coordinate
directions takes $d+1$ evaluations -- one at $\vs_i$ and one for each perturbation. Thus, approximating $\mJ$
from \eqref{eq:jacmat} takes $k(d+1)$ function evaluations. The potential benefits of revealing the directions of
variability may justify this cost, particularly if one is faced with a number of evaluations of $f$ that is exponential
in $d$ to construct an accurate surrogate. 

If evaluations of $f$ are very expensive, then we want to obtain the eigenvectors of $\hat{\mC}$ with as few as
possible. We can potentially use fewer than $k(d+1)$ evaluations by employing recently developed methods for matrix
completion~\cite{Candes10} under the assumption that $\mJ$ is row rank deficient, which is equivalent to a rank
deficient $\hat{\mC}$; see \eqref{eq:ctilde}. If $\mathrm{rank}\,(\mJ)=a$ -- corresponding to a function $f$ with $a$
directions of variability -- then we can recover $\mJ$ to within the precision of the finite difference approximation
by computing a constant times $a(k+d)$ entries of $\mJ$.

Let $\mathcal{D}$ be a subset of the pairs of indices $(i,j)$ with $i=1,\dots,d$ and $j=1,\dots,k$, where $i$
indexes the coordinates $s_i$ and $j$ indexes the design points $\vs_j$ from \eqref{eq:jacmat}. For a
$d\times k$ matrix $\mX$, define $P_\mathcal{D}(\mX)$ to return a vector of the entries of $\mX$ corresponding to the
index pairs in $\mathcal{D}$. For a given tolerance $\varepsilon$, the singular value thresholding (SVT)
algorithm~\cite{Cai10} seeks a solution to the convex optimization problem
\begin{equation}
\label{eq:matcomp}
\begin{array}{ll}
\mini{a} & \|\mX\|_\ast \\
\st       & \|P_\mathcal{D}(\mX)-P_\mathcal{D}(\mJ)\|\leq \varepsilon,
\end{array}
\end{equation}
where $\|\cdot\|_\ast$ is the nuclear norm. Note that the finite difference parameter for the approximate Jacobian
provides a natural tolerance on the constraints of the convex optimization problem. 

In fact, the SVT algorithm returns approximate singular vectors/values for $\mJ$, which saves the trouble of forming
$\hat{\mC}$ with $\mJ$; see \eqref{eq:ctilde}. The left singular vectors of $\mJ$ approximate the eigenvectors $\mV$,
and the singular values approximate the square roots of the eigenvalues; see \eqref{eq:eigdecomp}. We can use the
output of the SVT method directly to obtain the directions of variability. We demonstrate this approach in the
numerical examples in section \ref{sec:examples}. 

We will not comment on the cost of the SVT algorithm; we mention it for the case when computing more entries of $\mJ$
through evaluations of $f$ is more expensive than running the SVT algorithm. For example, if $f$ is evaluated with an
expensive PDE simulation, and the dimensions of $\mJ$ are in the tens to thousands, then this approach is appropriate. 

\subsection{Sampling from $\Omega_a$}
\label{sec:linprog}
To sample from the reduced space $\Omega_a$ defined in \eqref{eq:reddom}, we use a simple acceptance/rejection scheme.
We first determine an $a$-dimensional hyperrectangle that contains $\Omega_a$ by solving $a$ independent linear
programs,
\begin{equation}
\label{eq:hypprog}
\begin{array}{cc}
\mini{\vs} & \vv_i^T\vs, \\
\st & \vs_l\leq \vs \leq \vs_u,
\end{array}
\end{equation}
where $\vv_i$ is the $i$th column of $\mV_a$. Let $\vs_i^\ast$ be the minimizer of \eqref{eq:hypprog}. Then we define
the hyperrectangle $\tilde{\Omega}_a$ as
\begin{equation}
\tilde{\Omega}_a = \left\{\vt=\bmat{t_1 \\ \vdots \\ t_a} \;:\; 
\bmat{\vv_1^T\vs_1^\ast \\ \vdots \\ \vv_a^T\vs_a^\ast} \leq
\bmat{t_1 \\ \vdots \\ t_a} \leq
\bmat{-\vv_1^T\vs_1^\ast \\ \vdots \\-\vv_a^T\vs_a^\ast} \right\}.
\end{equation}
Notice that $\Omega_a\subset\tilde{\Omega}_a$, and we expect that the volume of the enclosing hyperrectangle will
be much larger than the volume of $\Omega_a$ in high dimensions. 

To draw a sample $\vs_a$ from $\Omega_a$, we draw $\vt$ uniformly from $\tilde{\Omega}_a$. If $\mV_a\vt\in\Omega$, then
we set $\vs_a=\vt$. If $\mV_a\vt\not\in\Omega$, but there exists a $\vz$ such that $\mV_a\vt + \mV_b\vz\in\Omega$, then
we also set $\vs_a=\vt$. To determine if such a point exists, we can attempt to solve the linear program,
\begin{equation}
\label{eq:linfeas}
\begin{array}{cc}
\mini{\vs} & \mathbf{0}^T\vs, \\
\st & \vs_a=\mV_a^T\vs \\
 & \vs_l\leq \vs \leq \vs_u.
\end{array}
\end{equation}
If a point $\vs^\ast$ is found that satisfies the constraints, then $\vz=\mV_b^T\vs^\ast$. If such a $\vz$ does not
exist, then we reject $\vt$. Notice that the objective function in \eqref{eq:linfeas} is essentially meaningless; it is
merely used to set the problem in terms easily entered into a linear program solver. 

Each sample from $\Omega_a$ is used to evaluate $g(\vs_a)$ as in \eqref{eq:gfun}, which we use to construct a surrogate
on the low dimensional subspace. 

\subsection{A practical algorithm}
We have now discussed all the pieces in the procedure for approximating $f$ on the low dimensional manifold. 
\begin{enumerate}
  \item \textbf{Compute the directions.} If one can evaluate $\vj(\vs)$, choose points $\vs_i\in\Omega$
  with $i=1,\dots,k$ and compute
  \begin{equation}
  \mJ=\bmat{\vj(\vs_1) & \cdots & \vj(\vs_k)}, \qquad \frac{|\Omega|}{k}\mJ\mJ^T\;=\;\hat{\mC}\;=\;\mV\Lambda\mV^T.
  \end{equation}
  If one can only evaluate $f(\vs)$, use the procedure from section \ref{sec:apjac} to approximate the
  eigendecomposition of $\hat{\mC}$.
  \item \textbf{Determine the directions of variability.} Examine the eigenvalues $\lambda_1,\dots,\lambda_d$ and choose
  a truncation $a<d$ according to their magnitude. (This judgment can be difficult to make algorithmically.) Set
  $\mV_a$ to be the first $a$ eigenvectors.
  \item \textbf{Evaluate $g$ at points in the reduced domain.} Use the acceptance/rejection method from section
  \ref{sec:linprog} to choose a set of design points $\vy_j\in\Omega_a$
  for $j=1,\dots,n$. For each design point, compute $g_j=g(\vy_j)$ using \eqref{eq:gfun}. Note that each evaluation may
  require the computation of $\vz_j$ using the approach described in section \ref{sec:linprog}. 
  \item \textbf{Approximate $f$ at a point in $\Omega$.} For a point $\vs\in\Omega$, compute $\vs_a=\mV_a^T\vs$.
  Approximate $g(\vs_a)$ using an interpolation procedure on the points $\{\vy_j\}$ and evaluations $\{g_j\}$. This
  approximation occurs on the space $\Omega_a$ of reduced dimension. Set $f(\vs)$ to be the approximation of $g(\vs_a)$. 
\end{enumerate}
Once the first three steps have been completed, the last step can be repeated as needed. For example, numerical
integration or optimization can be performed on $f(\vs)$ using the surrogate constructed on the reduced space
$\Omega_a$.

\section{Numerical Examples}
\label{sec:examples}
In this numerical exercise, we perform an uncertainty study on an elliptic PDE with a
random field model for the coefficients. Such problems are common test cases for methods in uncertainty
quantification~\cite{Babuska01,Ghanem91}.

\subsection{PDE model, input parameters, and quantity of interest}
\begin{figure}
\begin{center}
\subfloat[]{
\includegraphics[scale=0.25]{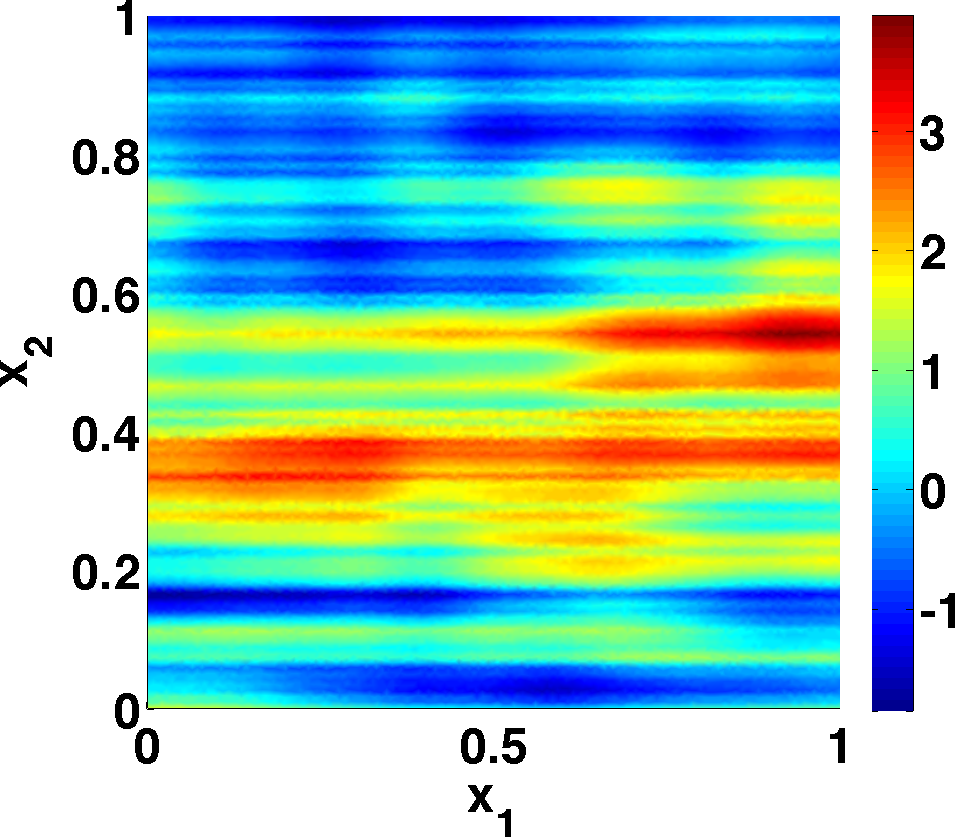}
}
\subfloat[]{
\includegraphics[scale=0.25]{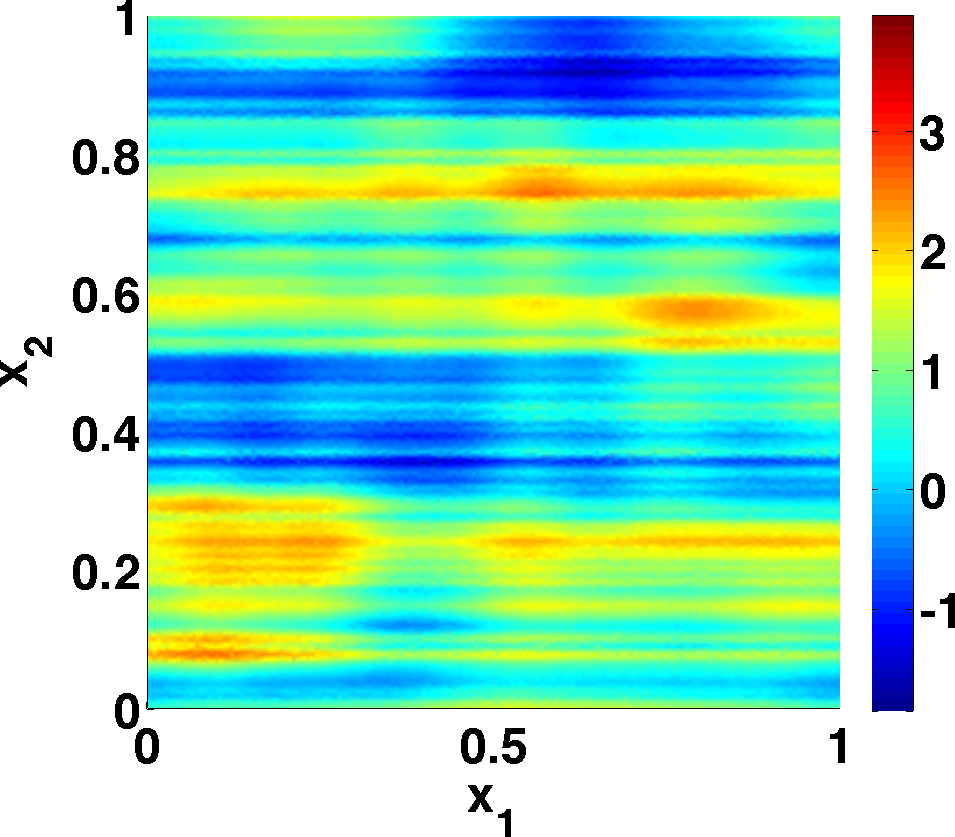}
}
\subfloat[]{
\includegraphics[scale=0.25]{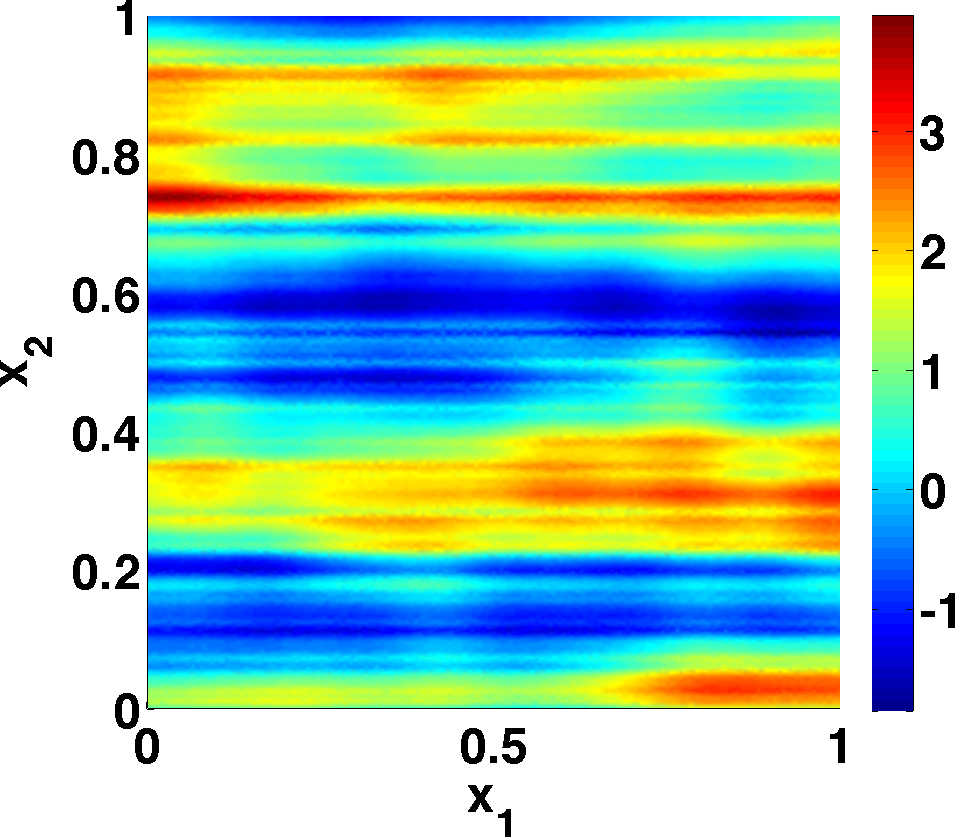}
}
\end{center}
\caption{Three realizations of the coefficients $\log(\alpha(\vx,\vs))$.}
\label{fig:a}
\end{figure}

\begin{figure}
\begin{center}
\subfloat[]{
\includegraphics[scale=0.25]{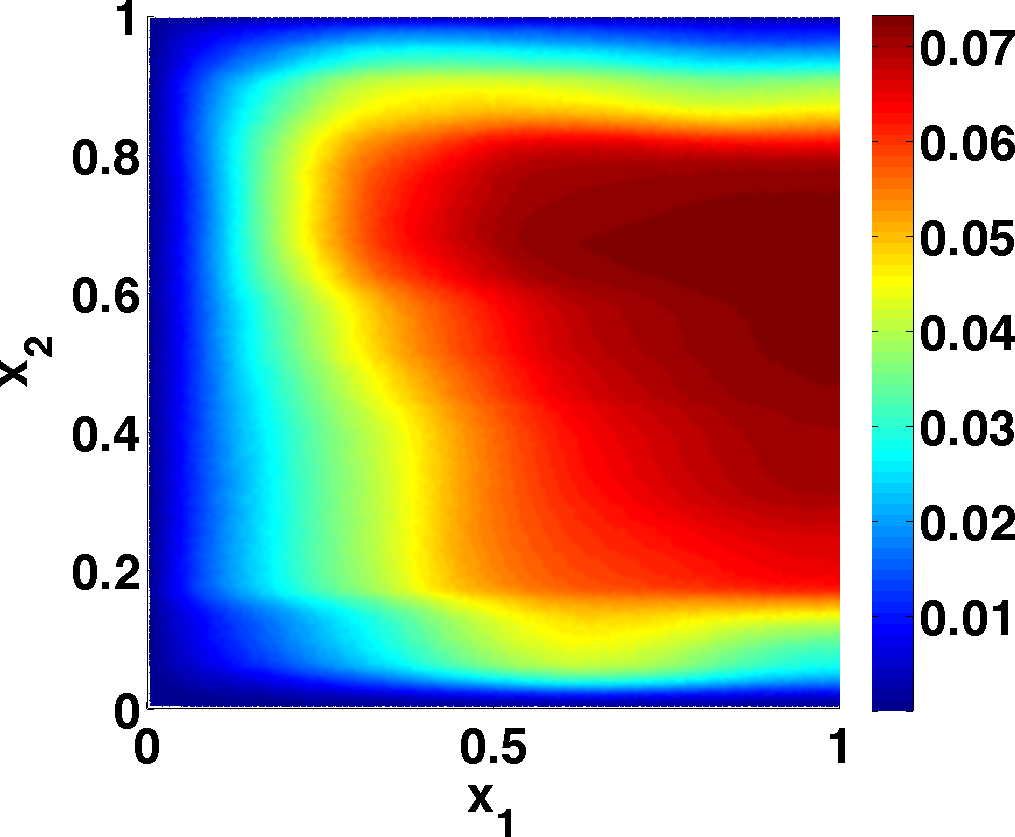}
}
\subfloat[]{
\includegraphics[scale=0.25]{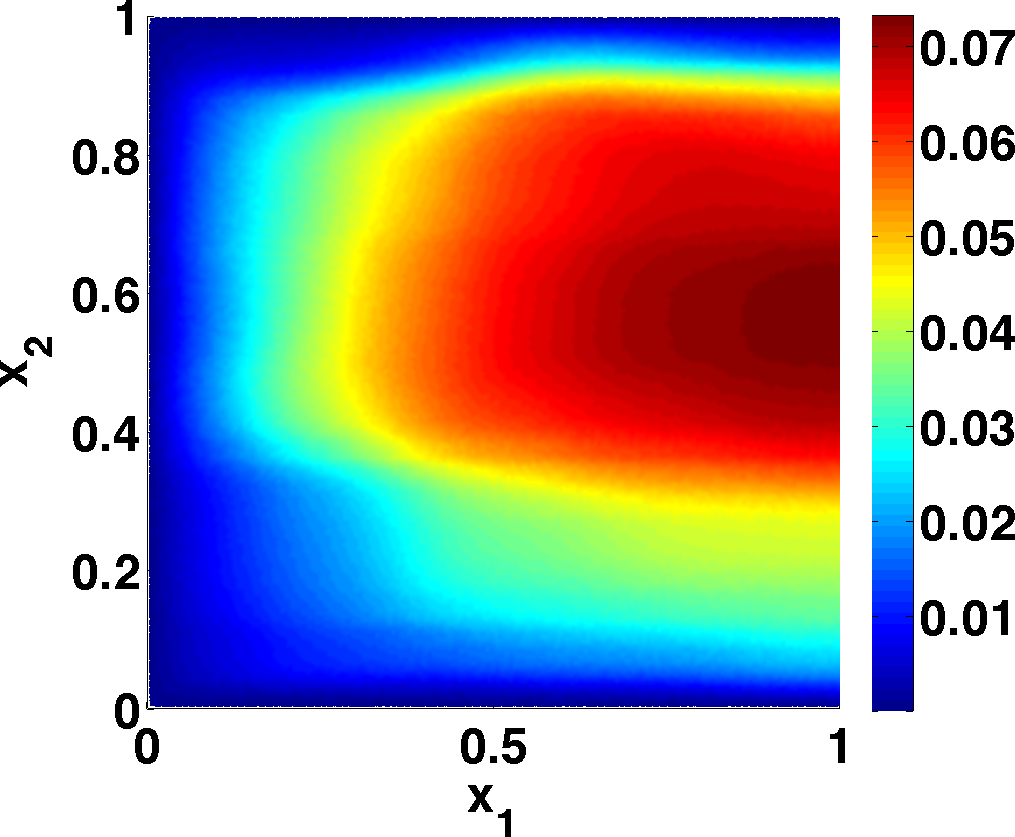}
}
\subfloat[]{
\includegraphics[scale=0.25]{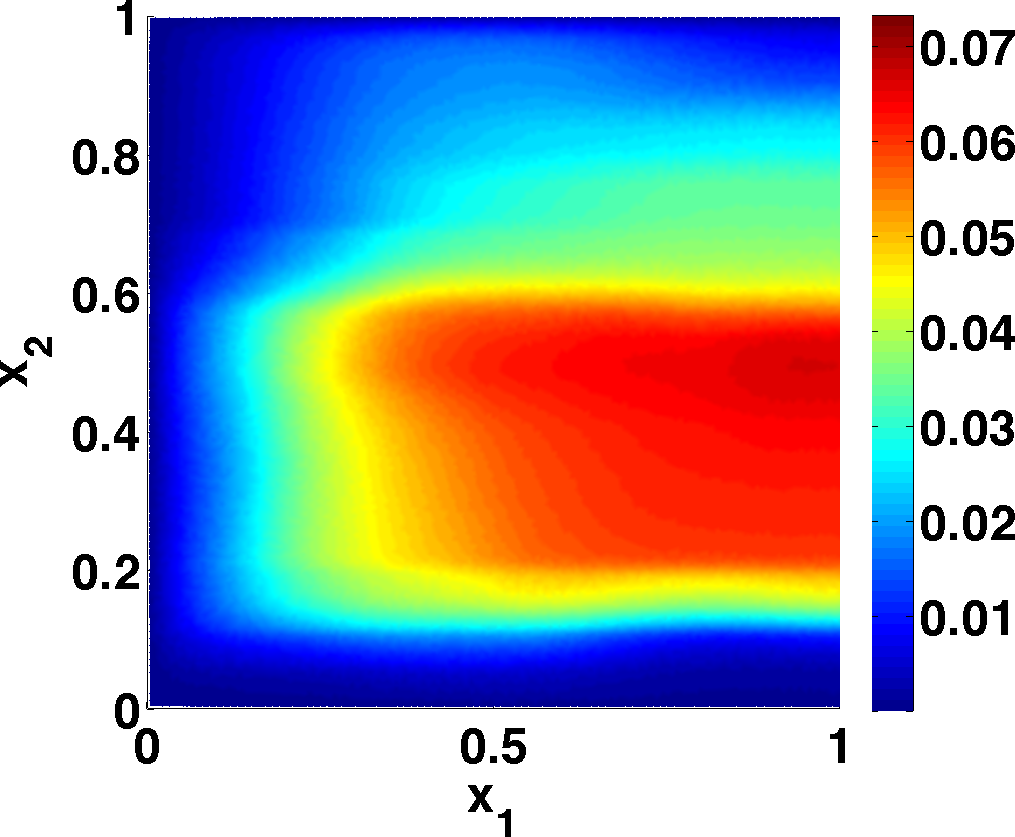}
}
\end{center}
\caption{Three realizations of the solution $u(\vx,\vs)$.}
\label{fig:u}
\end{figure}

Consider the following linear elliptic PDE. Let $u=u(\vx,\vs)$ satisfy
\begin{equation}
-\nabla\cdot(\alpha\nabla u) =1
\end{equation}
on the spatial domain $\vx\in[0,1]^2$. We set homogeneous Dirichlet boundary conditions on the left, top, and bottom of
the domain; denote this boundary by $\Gamma_1$. The right side of the domain -- denoted $\Gamma_2$ -- has a homogeneous
Neuman boundary condition. The log of the coefficients $\alpha=\alpha(\vx,\vs)$ of the differential operator are given
by a truncated Karhunen-Loeve type expansion
\begin{equation}
\label{eq:kl}
\log(\alpha) = \sum_{i=1}^d \phi_i\sqrt{\sigma_i} s_i,
\end{equation}
where the $s_i$ are independent, identically distributed uniform random variables on $[-2,2]$, and the 
$\{\phi_i,\sigma_i\}$ are the eigenpairs of the covariance operator
\begin{equation}
\mathcal{C}(\vx,\vy) = \exp\left\{-\left(\frac{(x_1-y_1)^2}{\rho_1} + \frac{(x_2-y_2)^2}{\rho_2} \right)\right\}
\end{equation}
with $\rho_1=1$ and $\rho_2=0.05$. The small $\rho_2$ models a short correlation length in the vertical coordinate. The
decay of the $\sigma_i$ justifies a truncation of $d=250$, so that the parameter space $\Omega$ for the problem is
the 250-dimensional hypercube $[-2,2]^d$. Three realizations of the log of the coefficients $\alpha$ and their
corresponding solutions $u$ are shown in figures \ref{fig:a} and \ref{fig:u}, respectively

Define the linear function $\tilde{Q}=\tilde{Q}(\vs)$ of the solution 
\begin{equation}
\tilde{Q} \;=\; \frac{1}{|\Gamma_2|} \int_{\Gamma_2} u \,d\vx.
\end{equation}
The quantity of interest for the uncertainty study is an approximate density function for $\tilde{Q}$.

\subsection{Finite element discretization}
Given a value for the input parameters $\vs$, we discretize the elliptic problem with a standard linear
finite element method using {\sc Matlab}'s PDE Toolbox. The discretized domain has 34320 triangles and 17361 nodes; the
eigenfunctions $\phi_i=\phi_i(\vx)$ from \eqref{eq:kl} are approximated on this mesh. The matrix equation for the
discrete solution $\vu=\vu(\vs)$ at the mesh nodes is
\begin{equation}
\label{eq:fem}
\mK\vu = \vf, 
\end{equation}
where $\mK=\mK(\vs)$ is symmetric and positive definite for all $\vs\in\Omega$. We can
approximate the linear functional $\tilde{Q}$ as
\begin{equation}
\label{eq:Q}
\tilde{Q} \;\approx\; \vc^T\vu \;=\; Q,
\end{equation}
where the elements of $\vc$ are zero except corresponding to nodes on $\Gamma_2$. The nonzero elements are constant
and scaled so that they sum to one; note that $\vc$ does not depend on $\vs$.  

\subsection{Adjoint variables for derivatives}
\begin{figure}
\begin{center}
\subfloat[]{
\includegraphics[scale=0.25]{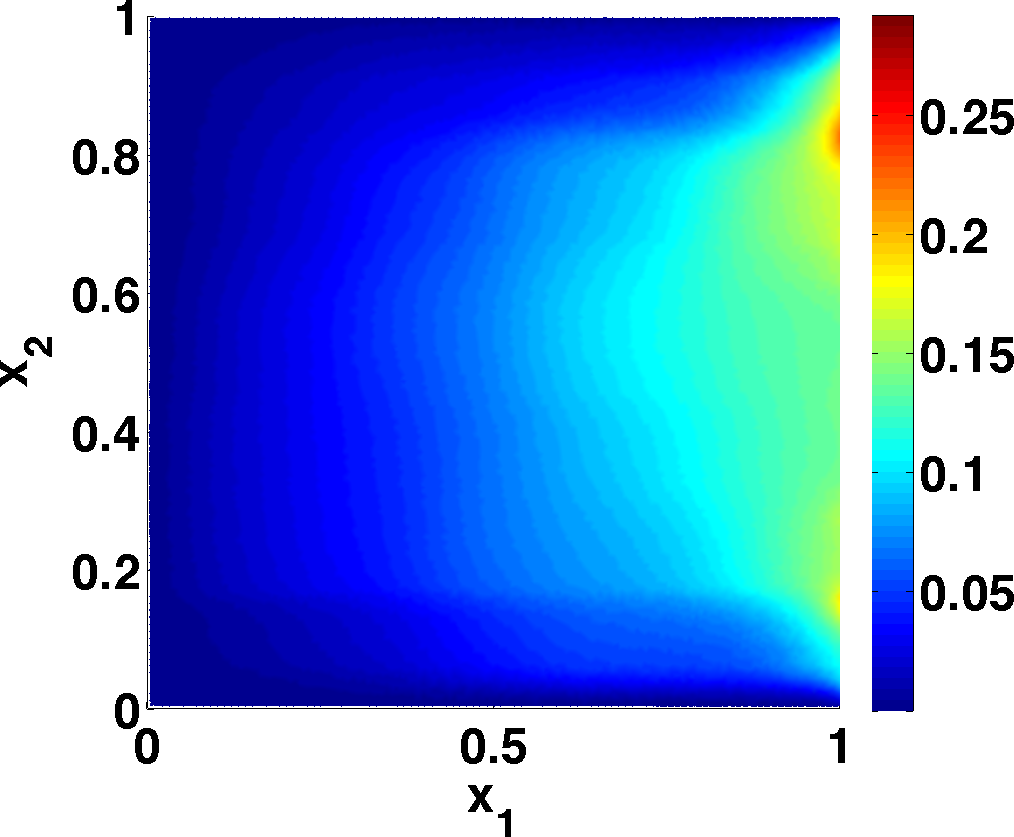}
}
\subfloat[]{
\includegraphics[scale=0.25]{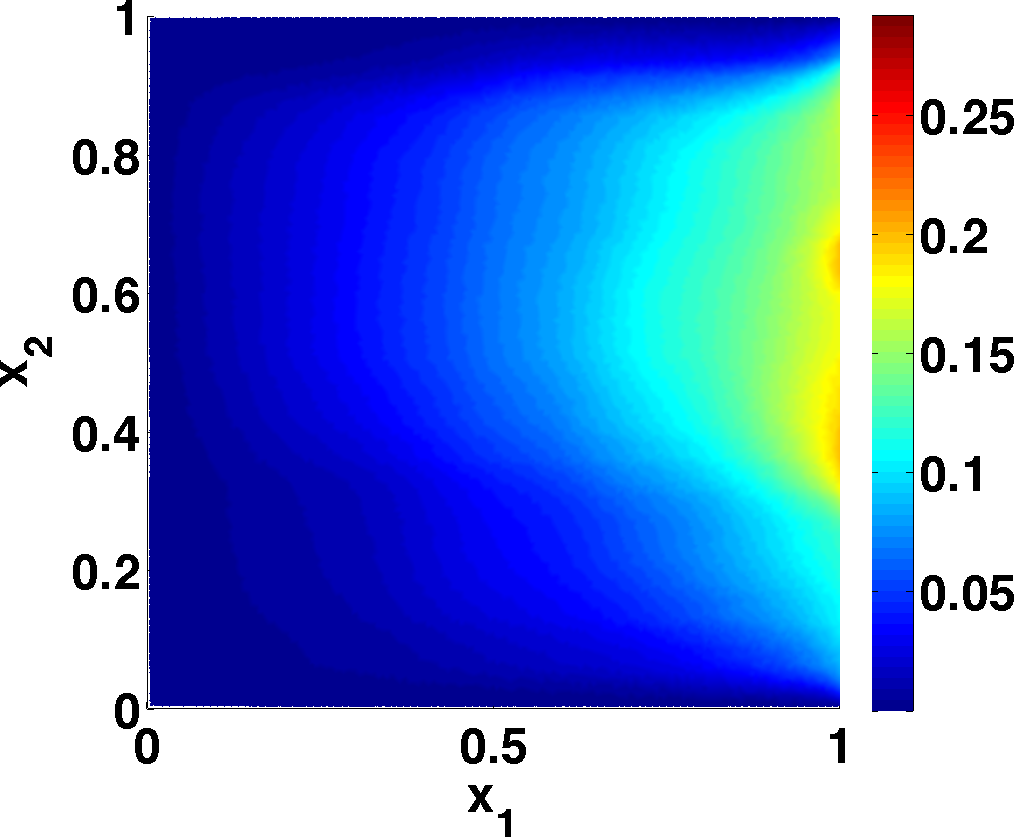}
}
\subfloat[]{
\includegraphics[scale=0.25]{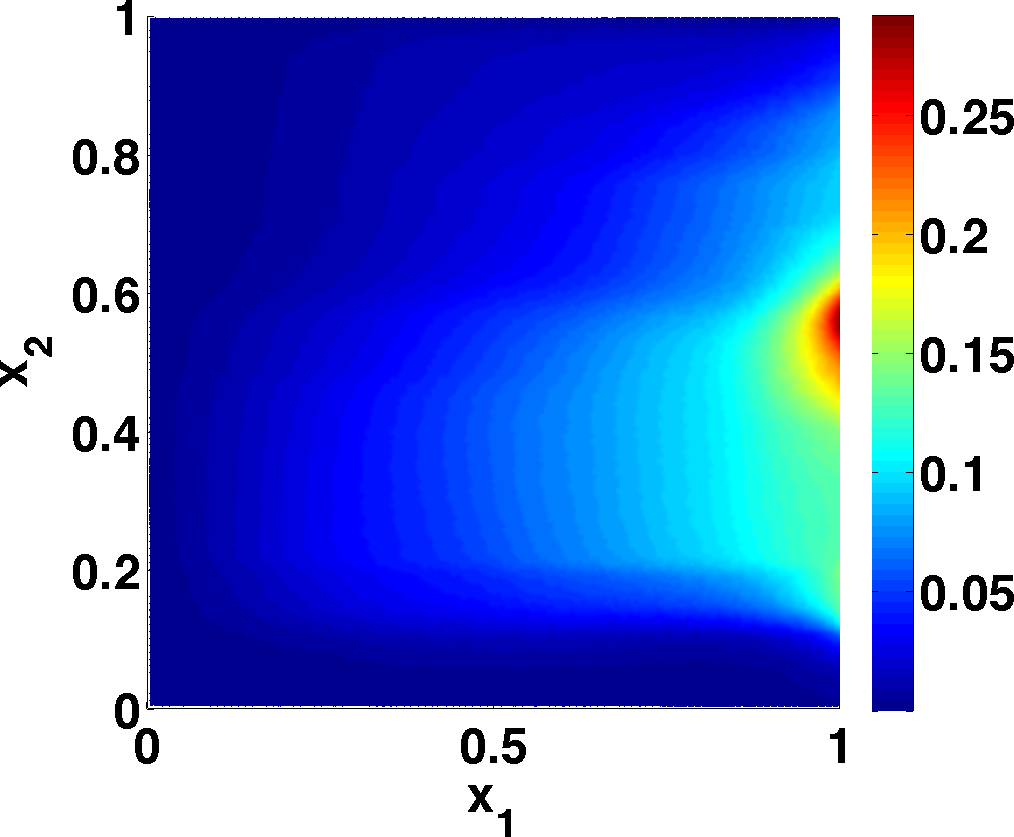}
}
\end{center}
\caption{Three realizations of the adjoint variables $\vy=\vy(\vs)$ on the mesh.}
\label{fig:y}
\end{figure}

Since the quantity of interest can be written as a linear functional of the solution, we can define adjoint variables
that we will help us compute the Jacobian of $Q$ with respect to the input parameters $\vs$. Notice that we can write
\begin{equation}
\label{eq:adjQ}
Q \;=\; \vc^T\vu \;=\; \vc^T\vu - \vy^T(\mK\vu - \vf),
\end{equation}
for any constant vector $\vy$. Taking the derivative of \eqref{eq:adjQ} with respect to the input $s_i$, we get
\begin{align*}
\frac{\partial Q}{\partial s_i} 
 &= \vc^T\left(\frac{\partial \vu}{\partial s_i}\right)
- \vy^T\left(\frac{\partial\mK}{\partial s_i}\vu + \mK\frac{\partial\vu}{\partial s_i}\right)\\
&= \left(\vc^T - \vy^T\mK\right)\left(\frac{\partial \vu}{\partial s_i}\right)
- \vy^T\left(\frac{\partial\mK}{\partial s_i}\right)\vu  
\end{align*}
If we choose $\vy$ to solve the adjoint equation 
\begin{equation}
\label{eq:adj}
\mK^T\vy=\vc, 
\end{equation}
then
\begin{equation}
\label{eq:derQ}
\frac{\partial Q}{\partial s_i} = - \vy^T\left(\frac{\partial\mK}{\partial s_i}\right)\vu.
\end{equation} 
Three realizations of the adjoint variables $\vy=\vy(\vs)$ are shown in figure \ref{fig:y}.

To approximate the Jacobian $\nabla Q$ at the point $\vs$, we compute the finite element solution with
\eqref{eq:fem}, solve the adjoint problem \eqref{eq:adj}, and compute the components with \eqref{eq:derQ}. The
derivative of $\mK$ with respect to $s_i$ is easy compute from the derivative of $a(\vx,\vs)$ and the same finite
element discretization. 

\subsection{Approximating the subspace}
\begin{figure}
\begin{center}
\includegraphics[scale=0.55]{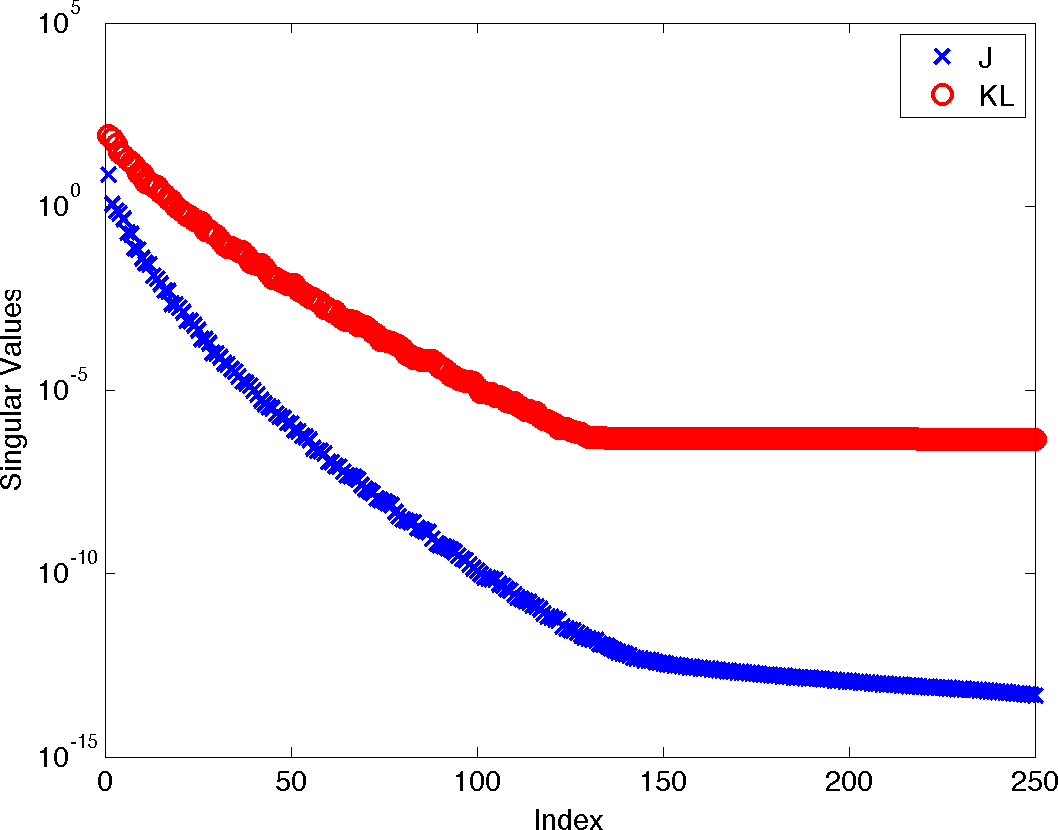}
\end{center}
\caption{The singular values of a reference computation of $\mJ$ with $10^4$ samples alongside the singular values of
the Karhunen-Loeve type expansion of the coefficients of the differential operator from \protect\eqref{eq:kl}.}
\label{fig:sigs}
\end{figure}

To apply the input reduction with the subspace detection technique, we first sample the Jacobian at random points in
$\Omega$ to construct $\mJ$ from \eqref{eq:jacmat}. From a reference computation of $10^4$ samples, we examine the
singular values to determine an appropriate truncation. The singular values of $\mJ$ are plotted in Figure
\ref{fig:sigs}; the decay justfies a truncation after five terms. For reference, we also plot the singular values
$\sqrt{\sigma_i}$ from the Karhunen-Loeve expansion \eqref{eq:kl}. The more rapid decay of the singular values of $\mJ$
shows that the particular output quantity of interest depends primarily on fewer variables than the correlated random
field modeling the coefficients of the differential operator.

In the remainder of the numerical exercise, we split the reference $10^4$ samples into five groups of 2000 samples. This
is to mimic an initial computational budget of 2000 samples, which we repeat five times to mildly alleviate affects
associated with a particularly good or bad choice of 2000 samples. Figures will display the results of each of the five
independent experiments.

\begin{figure}
\begin{center}
\subfloat[Reference error]{
\includegraphics[scale=0.32]{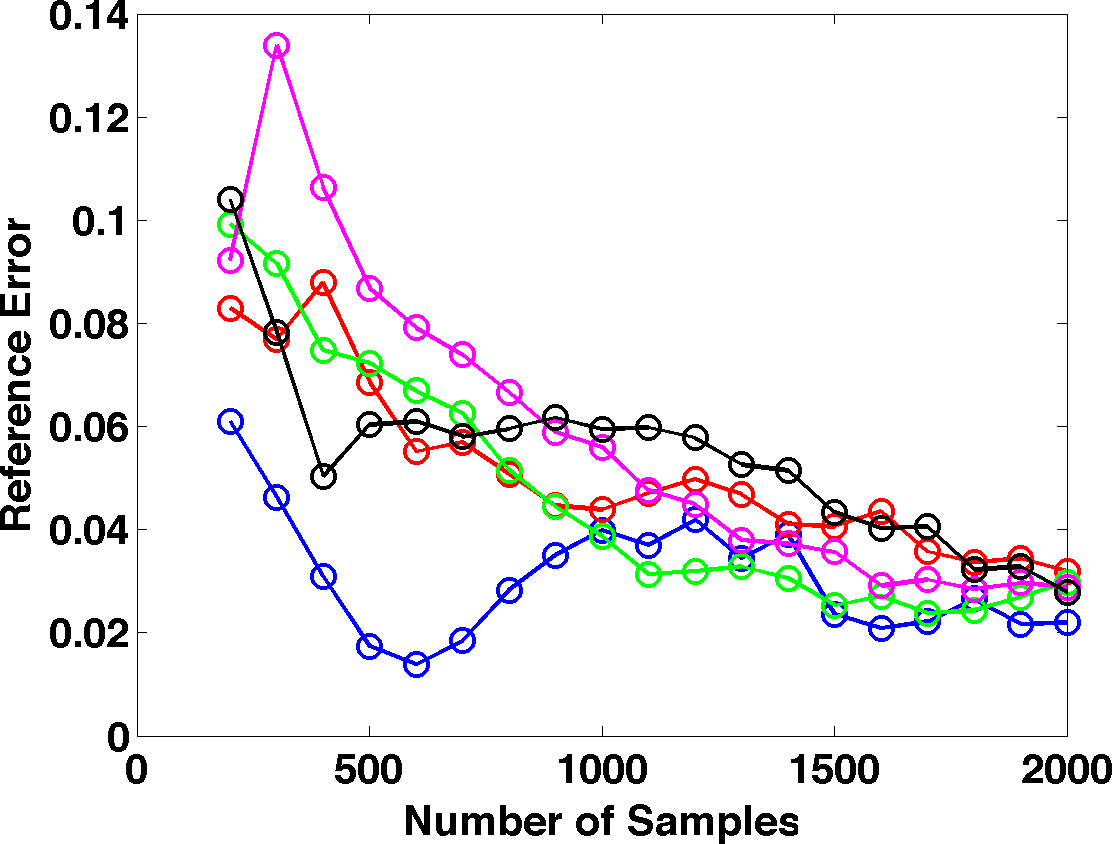}
\label{fig:abs}
}
\subfloat[Relative error]{
\includegraphics[scale=0.32]{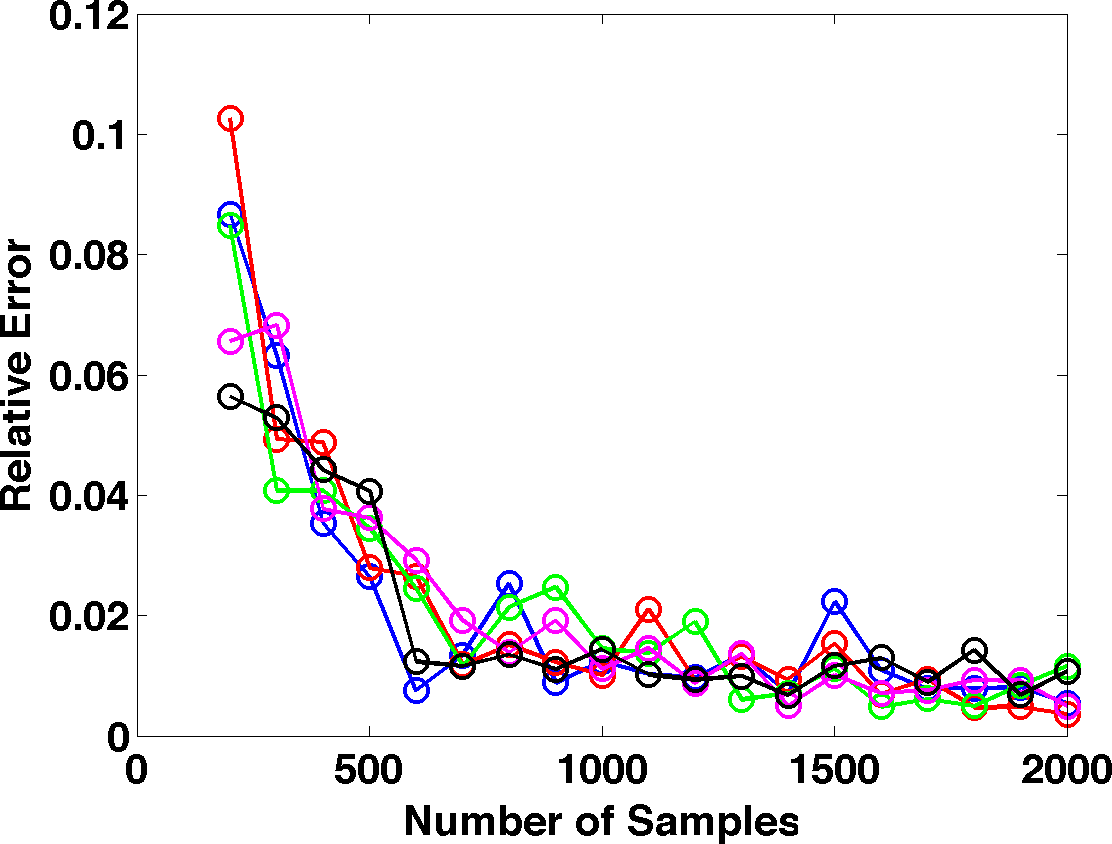}
\label{fig:rel}
}
\end{center}
\caption{The convergence of subspaces defined by taking the first $a=5$ left singular vectors from five different
random samplings of $\mJ$. The index $i$ corresponds to adding more samples to $\mJ$. Figure \protect\ref{fig:abs} shows
the error between the subspace from a reference solution of $10^4$ samples and the increasing the number of samples in
$\mJ$; see \protect\eqref{eq:abserr}. Figure \protect\ref{fig:rel} shows the relative change in the subspaces as more
samples are added; see \protect\eqref{eq:relerr}.}
\label{fig:subspaceerr}
\end{figure}

To check convergence of the projection onto the reduced subspace as more Jacobian samples are added, we compute the
left singular vectors of $\mJ$ for $m=100,200,300,\dots,2000$. The difference between the subspaces
defined by subsequent sets of samples $m_i$ and $m_{i+1}$ is given by
\begin{equation}
\label{eq:relerr}
\mathcal{E}_i^{\mathrm{rel}} = \|\mV_{a,m_{i+1}}\mV_{a,m_{i+1}}^T - \mV_{a,m_i}\mV_{a,m_i}^T\|,
\end{equation}
where $\mV_{a,m_i}$ are the first $a$ left singular vectors of $\mJ$ approximate with $m_i$ samples, and the norm is
the matrix 2-norm. In figures \ref{fig:abs}-\ref{fig:rel}, we plot $\mathcal{E}_i^{\mathrm{rel}}$ along with the error
of the projected subspace compared to the reference solution,
\begin{equation}
\label{eq:abserr}
\mathcal{E}_i^{\mathrm{abs}} = \|\mV_{a,m_i}\mV_{a,m_i}^T - \mV_{a,\mathrm{ref}}\mV_{a,\mathrm{ref}}^T\|,
\end{equation}
where $\mV_{a,\mathrm{ref}}$ are the first $a$ eigenvectors from the reference computation of $\mJ$. 

\begin{figure}
\begin{center}
\includegraphics[scale=0.55]{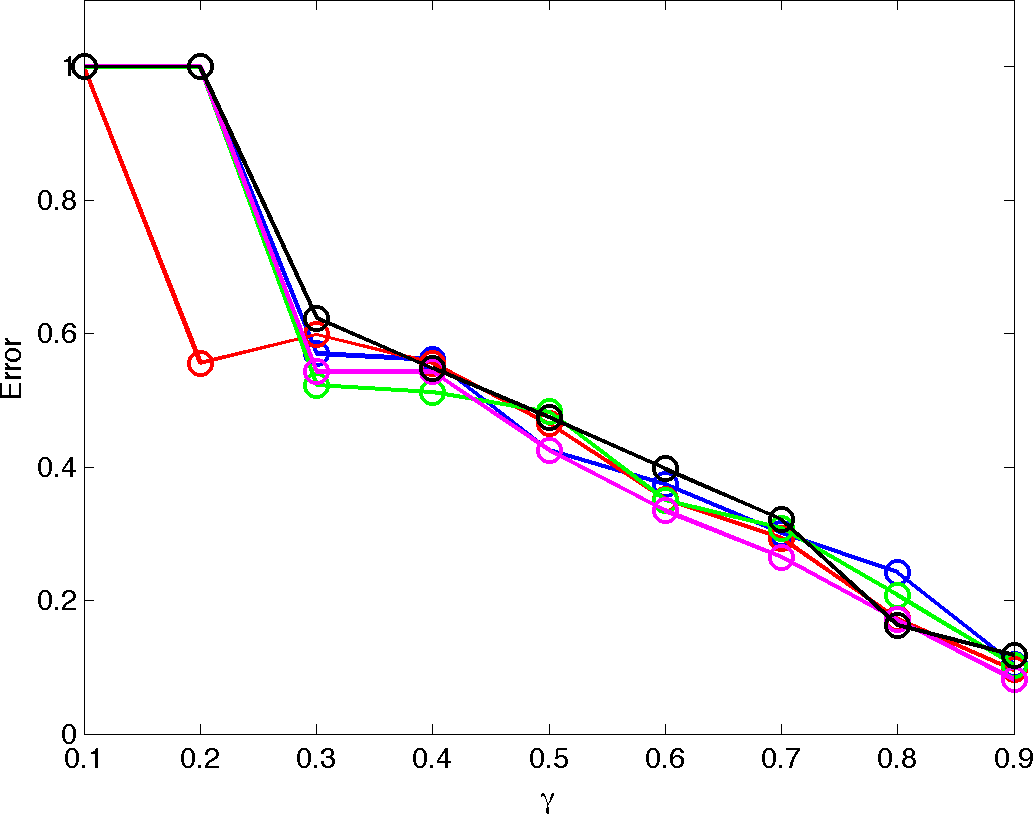}
\end{center}
\caption{The decrease in the error of the subspace defined by the first $a=5$ left singular vectors as more entries of
$\mJ$ with $m=2000$ samples are revealed according the proportion parameter $\gamma$; see \protect\eqref{eq:projerr}.}
\label{fig:projerr}
\end{figure}

Each column of $\mJ$ requires a forward solve, an adjoint solve, and 250 computations for the gradient $\nabla Q$. We
can try to reduce the number of derivative computations with the SVT algorithm for matrix completion, as described in
section \ref{sec:apjac}. Using the left singular vectors of $\mJ$ computed with $m=2000$ samples, we can test the SVT
method by uniformly subsampling the entries of $\mJ$. The number of subsampled entries is controlled by $\gamma$ with
$0<\gamma<1$, which is the proportion of entries revealed in the incomplete matrix. In figure \ref{fig:projerr}, we
plot the difference between the subspace from the subsampled $\mJ$ and the subspace from the full $\mJ$ for
$\gamma=0.1,0.2,\dots,0.9$,
\begin{equation}
\label{eq:projerr}
\mathcal{E}_\gamma = \|\mV_{a,\gamma}\mV_{a,\gamma}^T - \mV_{a}\mV_{a}^T\|.
\end{equation}
For the SVT algorithm, we used the {\sc Matlab} implementation from~\cite{Cai10} with the following parameters:
objective parameter \texttt{tau}=100, stopping criterion \texttt{tol}=1e-4, noise constraint \texttt{EPS}=1e-6, step
size \texttt{delta}=1, and maximum iterations \texttt{maxiter}=1000.

\subsection{Building the low dimensional surrogate}
We first map the 2000 design sites $\vs_i\in\Omega$ to get an initial set of design sites
$\vy_i=\mV_a^T\vs_i\in\Omega_a$; we have the evaluations of $Q$ associated with these points from the initial sampling
of the Jacobian. Following the process outlined in section \ref{sec:linprog}, we sample uniformly from the space
$\Omega_a$ to find 5000 additional design sites in the lower dimension subspace. 
The acceptance rate of the acceptance/rejection scheme is roughly 35\% (averaged over the five identical experiments),
which validates the intuition about the small volume of $\Omega_a$ relative to its enclosing hyperrectangle. It took an
average of 19201 linear programs to get the 5000 samples.

For each sample, we evaluate $Q$ from \eqref{eq:Q}. This gives us a total of 7000 points in the lower dimensional
subspace -- 2000 from the original Jacobian evaluations and 5000 from sampling on the reduced subspace $\Omega_a$ -- on
which to construct a surrogate. We use the kriging toolbox DACE~\cite{dace02} to build a surrogate on the lower
dimensional space $\Omega_a$.

\begin{figure}
\begin{center}
\subfloat[Exp. 1]{
\includegraphics[scale=0.25]{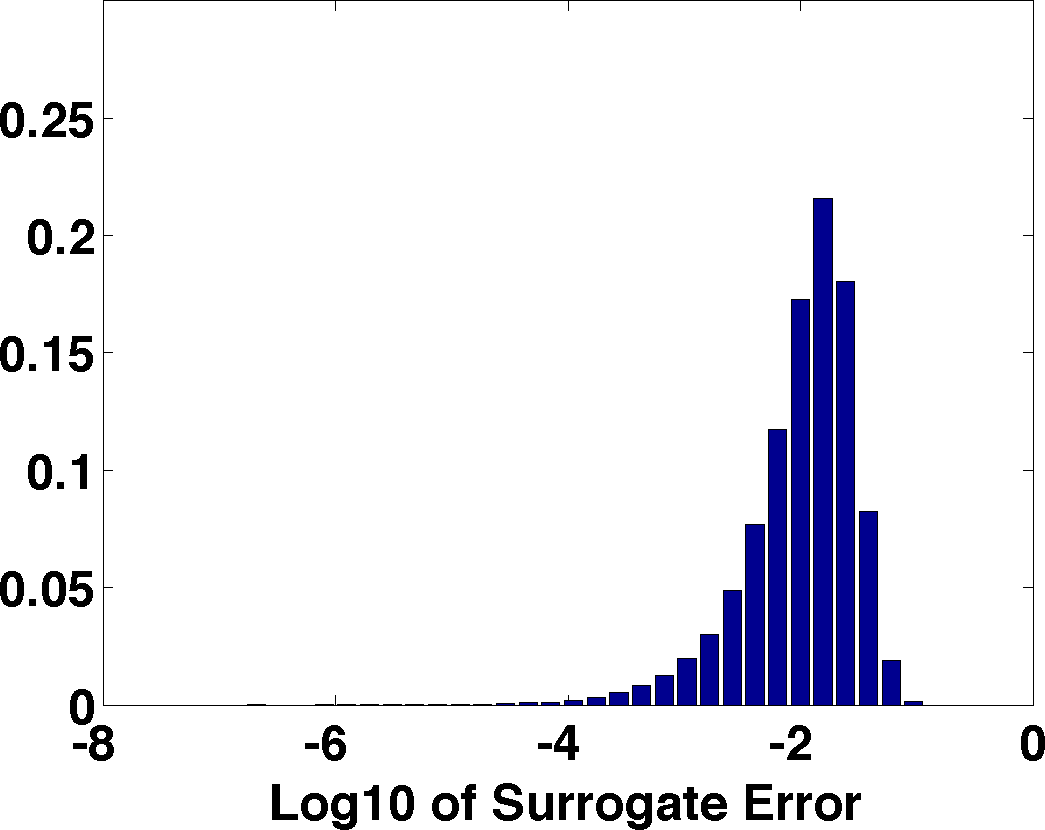}
}
\subfloat[Exp. 2]{
\includegraphics[scale=0.25]{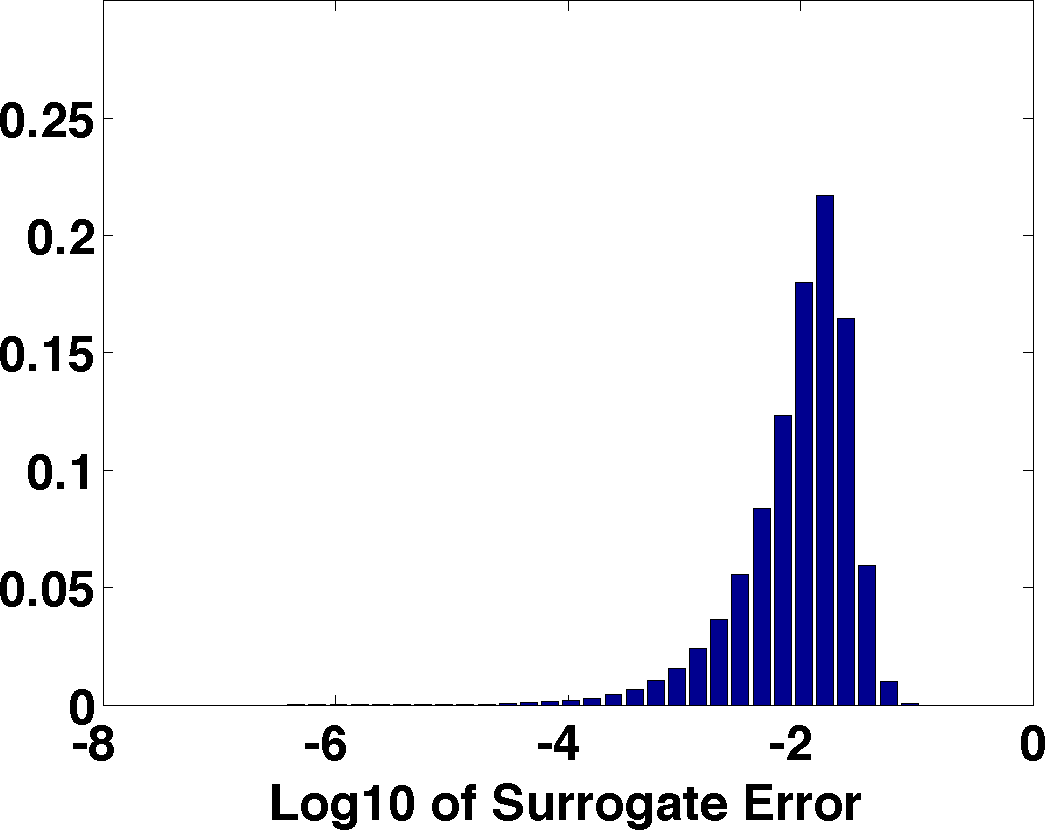}
}
\subfloat[Exp. 3]{
\includegraphics[scale=0.25]{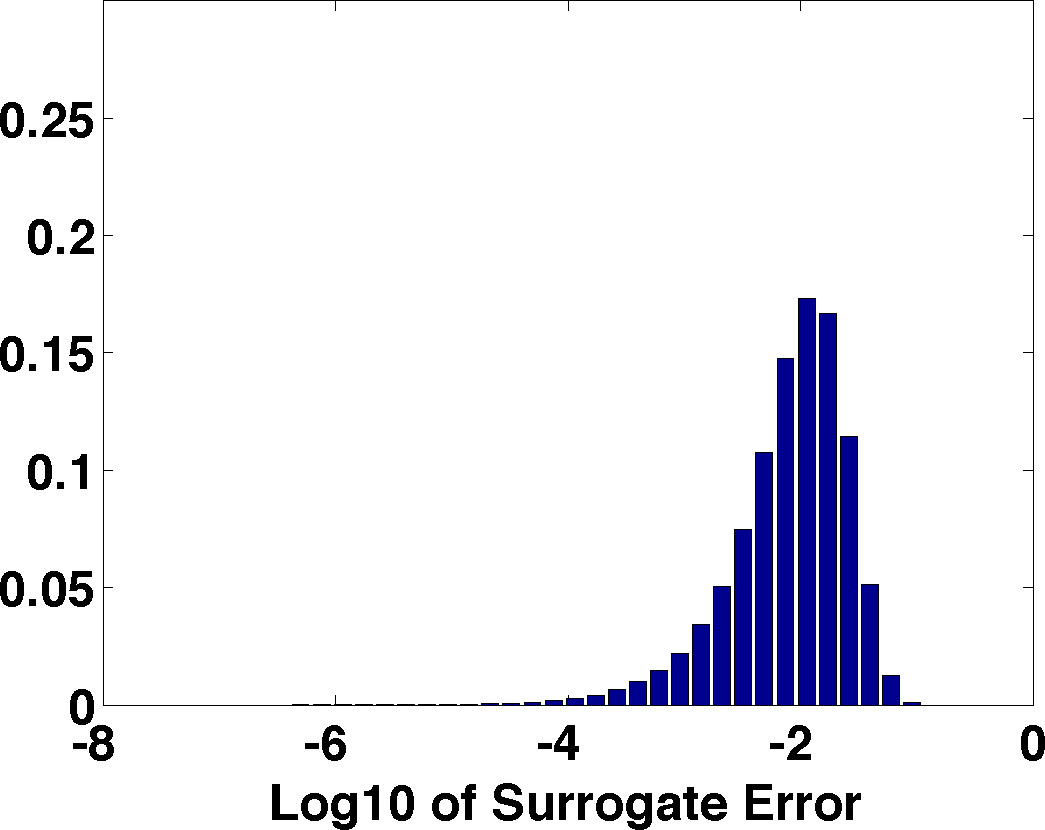}
}\\
\subfloat[Exp. 4]{
\includegraphics[scale=0.25]{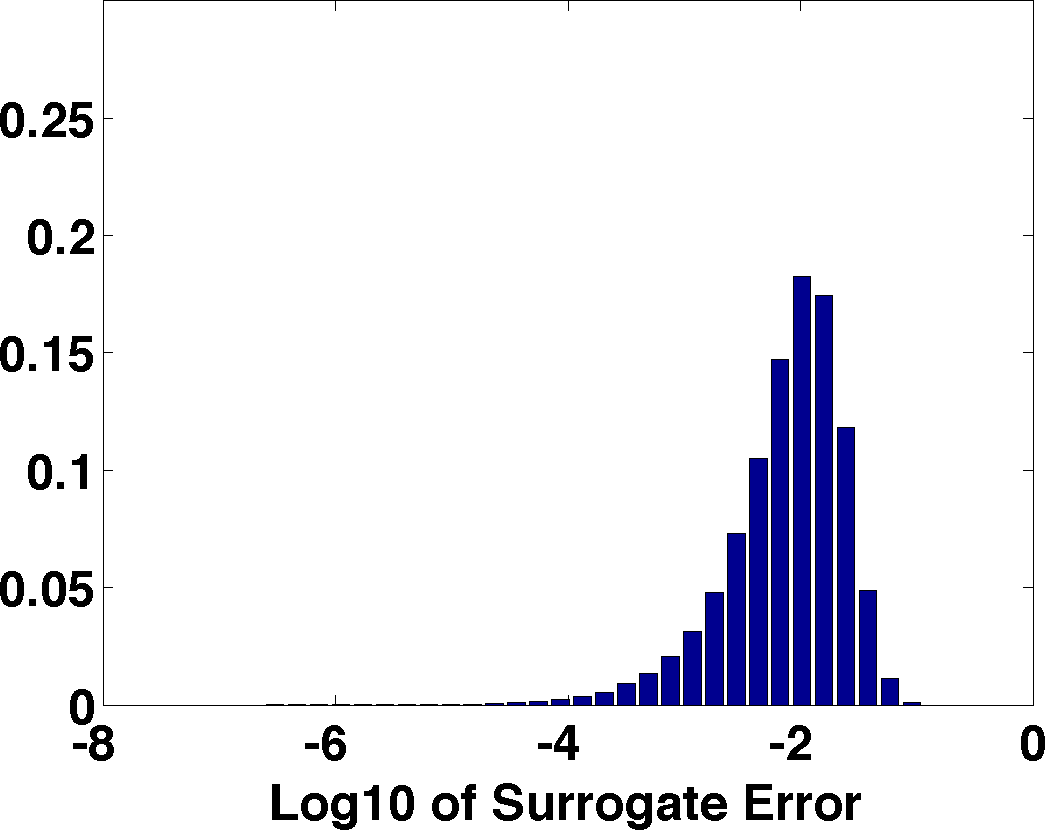}
}
\subfloat[Exp. 5]{
\includegraphics[scale=0.25]{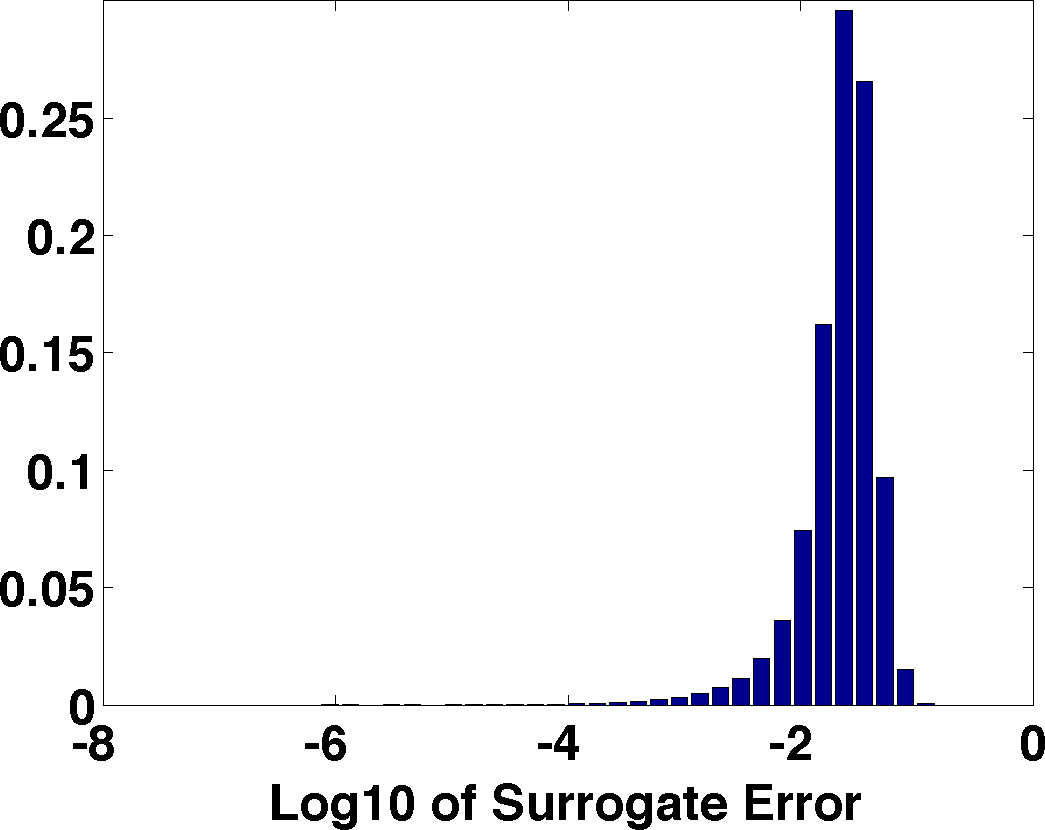}
}
\end{center}
\caption{Histograms of the log10 of the error between the surrogate approximation of $Q$ and the true value of $Q$ at
$10^5$ randomly sampled points from the full space $\Omega$. Each histogram corresponds to a different random sampling
used to construct the surrogate.}
\label{fig:histerr}
\end{figure}

Since the evaluation of $Q$ is relatively inexpensive in this example, we compare the surrogate's prediction of $Q$
with the actual $Q$ on $10^5$ points chosen uniformly at random from $\Omega$. The histograms of the log of the
surrogate error are shown in figure \ref{fig:histerr} -- one for each of the five experiments.  

\subsection{Approximating the desnity function}
\begin{figure}
\begin{center}
\subfloat[Full]{
\includegraphics[scale=0.25]{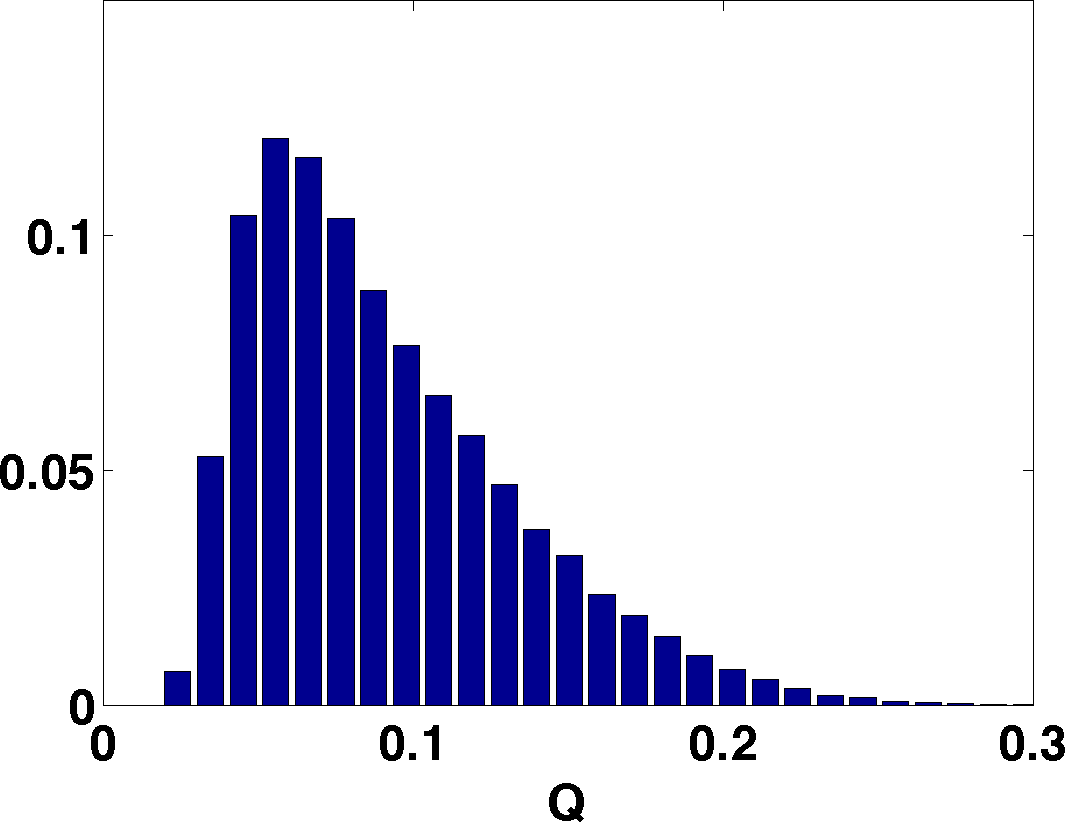}
\label{fig:truehist}
}
\subfloat[Exp. 1]{
\includegraphics[scale=0.25]{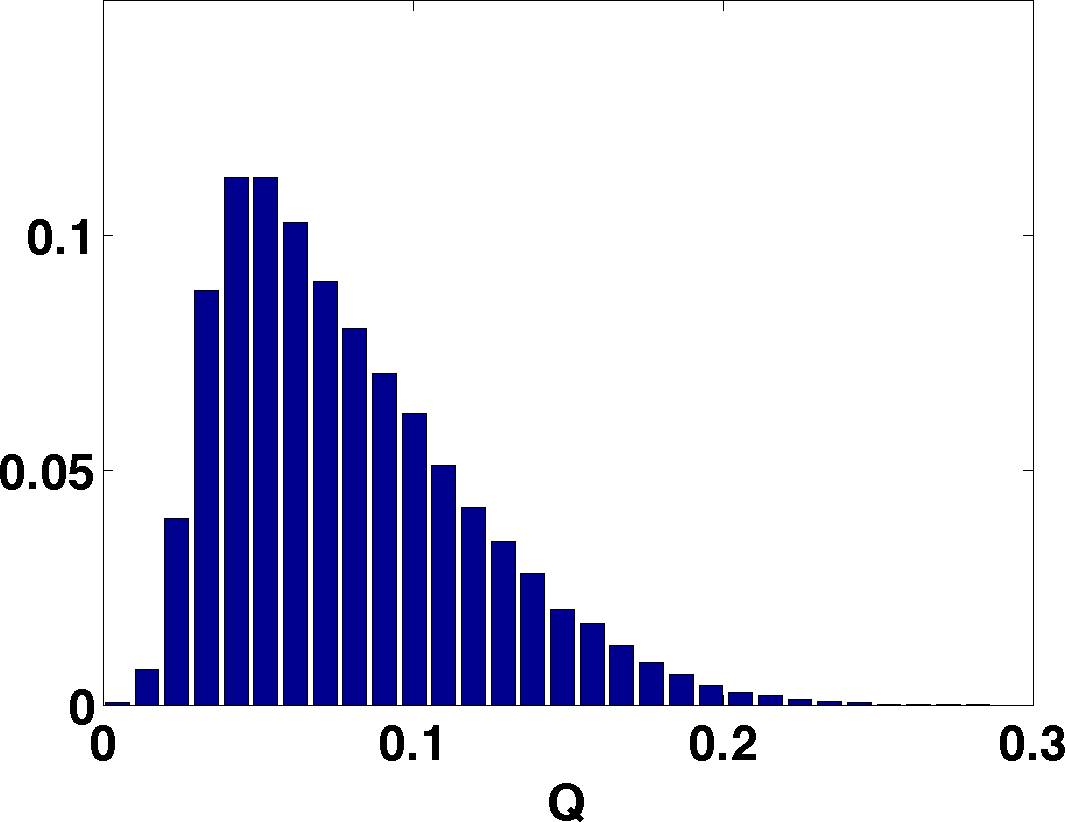}
}
\subfloat[Exp. 2]{
\includegraphics[scale=0.25]{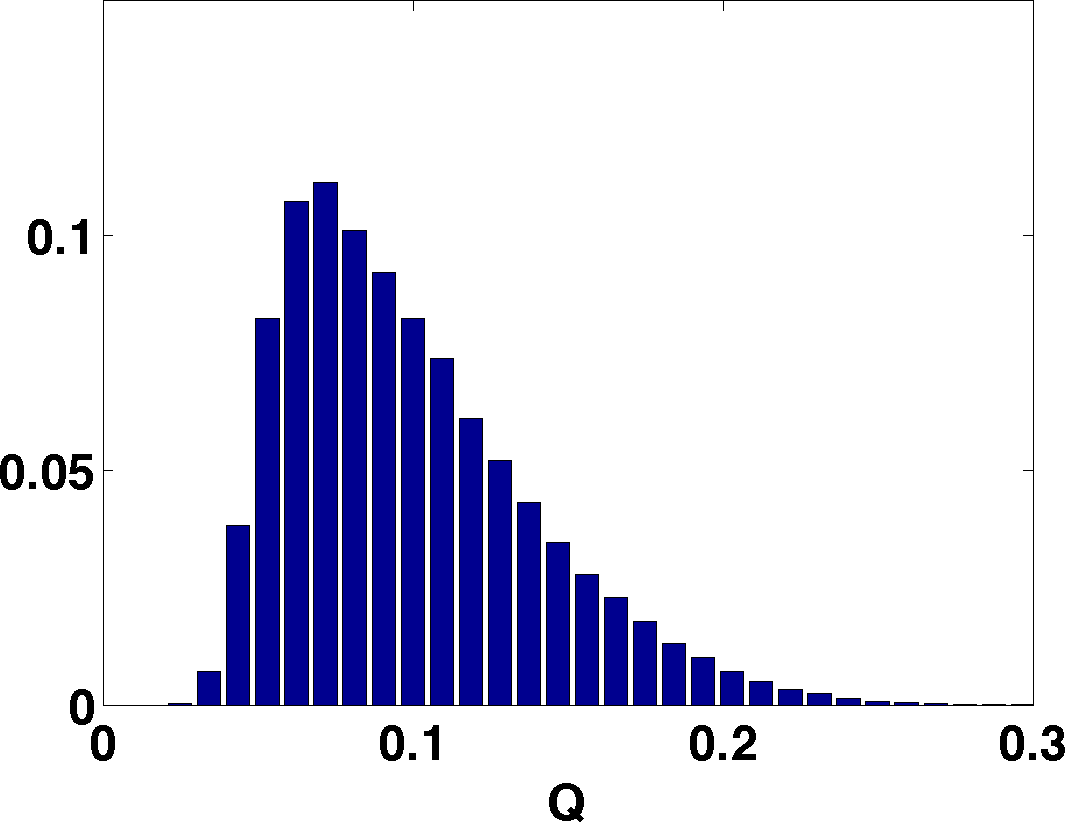}
}\\
\subfloat[Exp. 3]{
\includegraphics[scale=0.25]{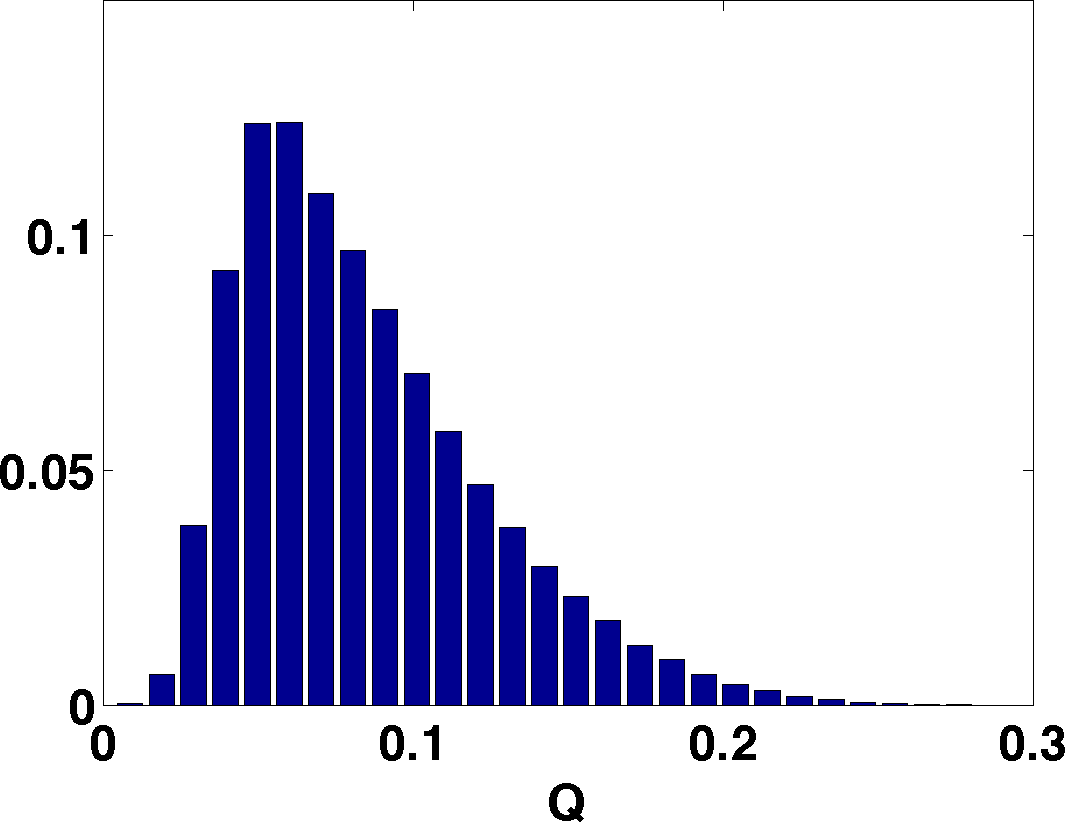}
}
\subfloat[Exp. 4]{
\includegraphics[scale=0.25]{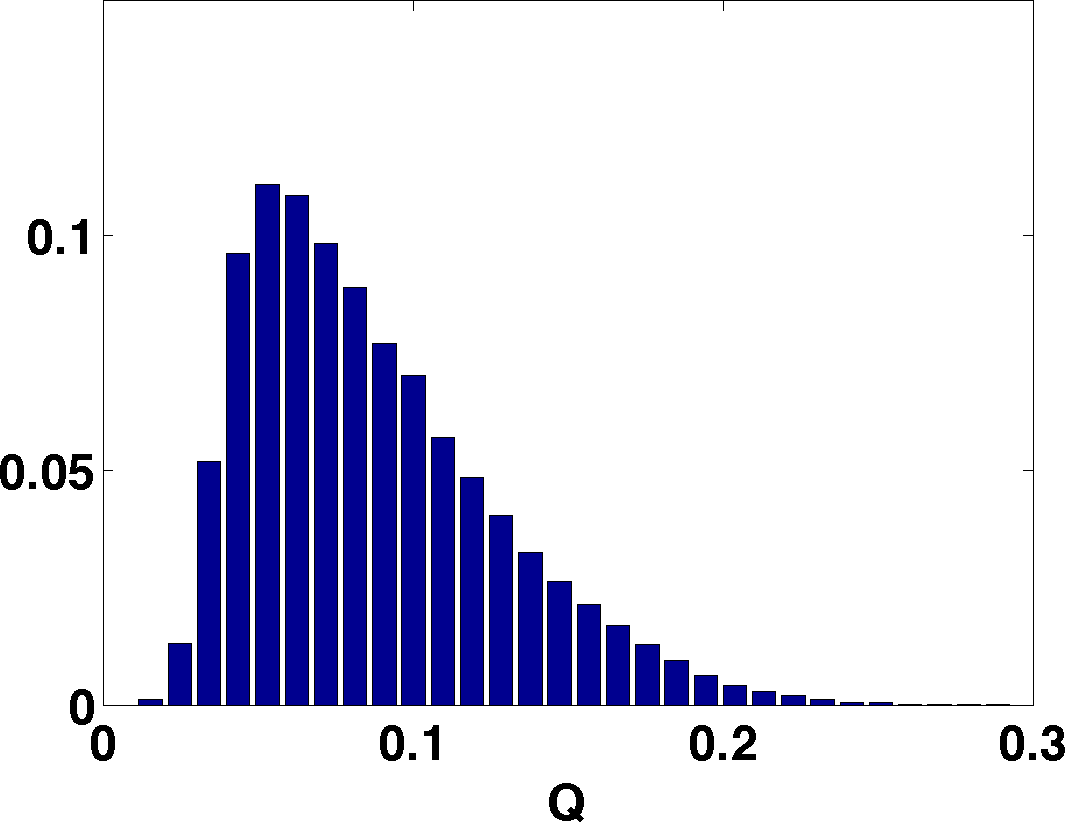}
}
\subfloat[Exp. 5]{
\includegraphics[scale=0.25]{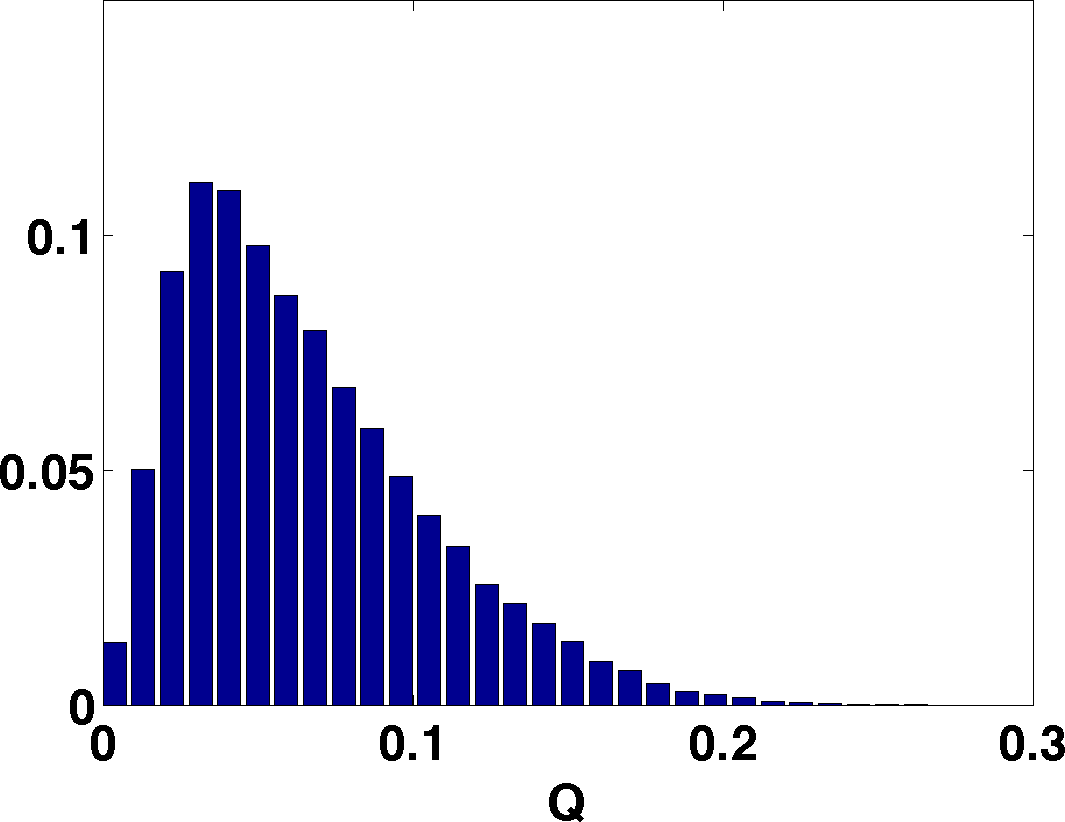}
}
\end{center}
\caption{Histograms of the quantity of interest $Q$ with $10^5$ samples from the input space $\Omega$. The top right
histogram \protect\ref{fig:truehist} is the histogram of the full model. The others are for the five experiments from
the low dimensional surrogate.}
\label{fig:hist}
\end{figure}

To approximate the density function of $Q$, we draw samples from the full space $\Omega$, use the low
dimensional surrogate to approximate the output quantity of interest, and build a histogram of the samples.
Specifically, for a point $\vs\in\Omega$, we compute $\vs_a=\mV_a^T\vs$. Then we use the low dimensional surrogate to
approximate $Q$ at $\vs_a$. Using $10^5$ such evaluations, we obtain reasonably well-converged histograms. In figure
\ref{fig:hist}, we plot the histogram of $10^5$ evaluations of $Q$ from the full model alongside the histograms from each
surrogate experiment. We see that surrogate approximates the bulk of the histogram reasonably but loses accuracy near
the tails. This is expected; the low dimensional subspace is detected by averaged variability. Therefore, we do not
expect to capture extremes of $Q$, and this is reflected in a loss of accuracy in the tails.

\section{Conclusion}

We have presented a method for detecting the primary directions of variability of a function of many variables. We have
described how to exploit these directions to construct a surrogate on a low dimensional subspace of the high dimensional
input space. We demonstrated this procedure on an uncertainty quantification study with a model problem of an elliptic
PDE with variable coefficients that depend on 250 independent input parameters.

\bibliographystyle{siam}
\bibliography{/home/paulcon/Dropbox/workspace/paulconstantine}

\end{document}